\documentclass[english]{article}

\usepackage{amsfonts}
\usepackage{amsmath}
\usepackage{graphicx}
\usepackage{epstopdf}
\usepackage{algorithmic}
\usepackage[latin1]{inputenc}
\usepackage{icomma}

\usepackage{hyperref}
\usepackage{cleveref} 
\usepackage{amsmath}
\usepackage{amsfonts,color}
\usepackage{mathtools}
\usepackage{tcolorbox}
\usepackage[pdf]{pstricks}
\usepackage{pst-node,pst-tree}

\usepackage{pstricks-add}
\usepackage[latin1]{inputenc}
\usepackage[bottom]{footmisc}
\usepackage{setspace}
\usepackage{boldline}
\usepackage{graphics}
\usepackage{graphicx}
\usepackage{url}
\usepackage{siunitx}
\usepackage{epsfig}
\usepackage{graphicx}
\usepackage{geometry}

\usepackage{authblk}

\usepackage{enumitem}
\usepackage{float}
\usepackage[all]{xy}
\usepackage{color}
\usepackage{multirow}
\usepackage{makecell}
\usepackage{hhline}
\usepackage{verbatim} 
\usepackage{geometry}

\definecolor{palegreen}{cmyk}{0.25,0,0.25,0}
\definecolor{lightgreen}{cmyk}{0.5,0,0.5,0}
\definecolor{darkgreen}{cmyk}{1,0,1,0.75}

\def\R{\mathbb{R}}
\newcommand{\RR}{{\mathbb R}} 
\newcommand{\NN}{{\mathbb N}} 

\def\susceptible{S}
\def\exposed{E}
\def\infected{I}
\def\infectedU{I^m}
\def\recovered{R}
\def\hospitalized{H}
\def\hospitalizedC{H^c}
\def\dead{D}

\def\duration{\gamma}
\def\fraction{\phi}
\def\rateg{\Lambda}
\def\state{x}
\def\stateh{\mathbf{x}}
\def\x{\mathbf{x}}
\def\u{\mathbf{u}}
\def\ic{\state_0}
\def\control{u}
\def\dynamics{f}
\def\prob{p}
\def\rate{\beta}

\def\Set{\mathcal{S}}
\def\feedback{\mathfrak{u}}
\def\feedbackf{\varphi}
\def\observation{\mathcal{O}}
\def\trigger{\mathfrak{m}}
\def\output{\mathbf{P}}
\def\outputc{P}
\def\tradeoff{\mathcal{C}}
\def\thresholdd{\mathcal{D}}
\DeclarePairedDelimiter\ceil{\lceil}{\rceil}
\newcommand{\subscript}[2]{$#1#2$}

\def\paso{\Delta t}

\newcommand{\outputcv}[2]{\mathbf{\hat P}_{#1:#2}}

\def\crref{c^{\rm ref}}

\def\it{t_0} 
\def\horizon{T} 
\def\hist{\Delta } 
\def\delay{t_{\rm min}} 
\newcommand{\titf}[2]{[\![{#1}\!:\!{#2}]\!]}
\def\ddfunction{\xi_{\rm day}} 

\def\expose2{E_2}

\def\infecte2{I^{severe}}
\def\recovered{R}

\def\indicator{\mathcal{I}}
\def\threshold{\theta}

\newcommand{\X}{\mathbf{X}}

\newcommand{\U}{\mathbf{U}}
\newcommand{\UU}{{\mathbb U}}

\usepackage{amsopn}

\newcommand{\XX}{{\mathbb X}}

\title{Assessment of event-triggered policies of nonpharmaceutical interventions based on epidemiological indicators}

\author[$\dagger$]{C. Castillo-Laborde}
\author[$\star$]{T. de Wolff}
\author[$\ddag$]{P. Gajardo\footnote{Corresponding author: pedro.gajardo@usm.cl}}
\author[$\ddag$]{R. Lecaros}
\author[$\triangleleft$]{G. Olivar}
\author[$\odot$]{H. Ramírez C.}

\affil[$\dagger$]{\footnotesize Centro de Epidemiología y Políticas de Salud (CEPS), Universidad del Desarrollo, Santiago, Chile}
\affil[$\star$]{\footnotesize Centro de Modelamiento Matemático (CNRS UMI 2807), Universidad de Chile, Santiago, Chile}
\affil[$\ddag$]{\footnotesize Departamento de Matemática, Universidad Técnica Federico Santa María, Valparaiso, Chile}

\affil[$\triangleleft$]{\footnotesize Departamento de Ciencias Naturales y Tecnología, Universidad de Aysén, Coyhaique, Chile}

\affil[$\odot$]{\footnotesize Departamento de Ingenier{\'i}a Matem{\'a}tica and Centro de Modelamiento Matem{\'a}tico (CNRS UMI 2807), Universidad de Chile, Santiago, Chile}

\date{\empty}

\begin{document}
\maketitle

\vspace{-.8cm}

\begin{abstract}
Nonpharmaceutical interventions (NPI) such as banning public events or instituting lockdowns have been widely applied around the world  to control the current COVID-19 pandemic. Typically, this type of intervention is imposed when an epidemiological indicator in a given population exceeds a certain threshold.  Then, the nonpharmaceutical intervention is lifted when the levels of the indicator used have decreased sufficiently.  What is the best indicator to use? In this paper, we propose a mathematical framework to try to answer this question. More specifically, the proposed framework permits to assess and compare different event-triggered controls based on epidemiological indicators. Our methodology consists of considering some outcomes that are consequences of the nonpharmaceutical interventions that a decision maker aims to make as low as possible.  The peak demand for intensive care units (ICU) and the total number of days in lockdown are examples of such outcomes. If an epidemiological indicator is used to trigger the interventions, there is naturally a trade-off between the outcomes that can be seen as a curve parameterized by the trigger threshold to be used. The computation of these curves for a group of indicators then allows the selection of the best indicator the curve of which   dominates the curves of the  other indicators. This methodology is illustrated using  indicators in the context of COVID-19 using deterministic compartmental models in discrete-time, although the framework can be  adapted for a larger class of models.\\

\end{abstract}

\section{Introduction}

After the initial outbreak of COVID-19 in Wuhan, China, the novel coronavirus SARS-COV-2 has spread to multiple countries \cite{Li2020} causing over 37 million cases and more than 1 million deaths by October 12, 2020 \cite{WorldHealthOrganizationWHO.SituationReport2020}. Since the emergence and global spread of COVID-19 in a context where no approved vaccines, treatments, or prophylactic therapies are available~\cite{Alvi2020}, countries have implemented different nonpharmaceutical interventions (NPI) \cite{Flaxman2020}. These prevention, containment, and mitigation policies have been implemented in order to flatten the peak of critical cases and consequently,  to prevent as far as possible the  health systems from being overwhelmed \cite{Alvi2020,OECD2020}. 
 
Public health policy responses to the pandemic aim both to limit the number and duration of social contact. They include measures such as early detection of cases, and tracing and isolating infected individual's contacts (testing and contact tracing), the use of Personal Protective Equipment (PPE) by health care workers, and social distancing measures, among others \cite{Alvi2020,OECD2020}. Social distancing interventions include working from home, school closures, cancellation of public gatherings and restrictions on unnecessary movement outside one's residence such as national or partial lockdowns \cite{Alvi2020,Davies2020,Flaxman2020,OECD2020}. These measures have  been shown to be successful in reducing the spread of SARS-COV-2 \cite{Alvi2020,Flaxman2020}. Lockdowns, in particular, have had a substantial effect on the transmission as measured by the changes in the estimated reproduction number  \cite{Flaxman2020}. In fact, a modeling study of the demand for hospital services in the UK found that intermittent periods of more intensive lockdown-type measures are predicted to be effective in preventing health system overload \cite{Davies2020}.

Beyond the aforementioned health consequences of the pandemic, it also has important economic consequences \cite{Martin2020,McKibbin2020}, throwing many countries into recession and possible economic depression \cite{Brodeur2020}. The economic impact can be direct and indirect \cite{Beland2020,Brodeur2020,Fan2018,Martin2020,Meltzer1999}, at the individual  or aggregate level, affecting firms and different sectors \cite{Martin2020,Nicola2020}, and with effects in the short, medium and long term  \cite{Beland2020,Brodeur2020,Gong2020,Nicola2020}. The negative economic effects may vary based on the severity of the social distancing measures, the length of their implementation, and the degree of compliance  \cite{Brodeur2020}, and, as usual, in times of crises the most vulnerable groups are affected the most \cite{Buheji2020,Marmot681}.

Health systems must face direct costs related to the diagnosis and confirmation testing, outpatient and inpatient care, and vaccines and pharmaceutical treatments when these become available  \cite{Fan2018,Meltzer1999}. On the other hand, productivity is affected by both absenteeism (patients and preventive isolation) and deaths caused by the virus \cite{Brodeur2020,Fan2018,Meltzer1999,UNDP2020}.

In this context, authorities are constantly searching for balance between health and economic consequences; with the latter being the result of not only the former but also of the measures or policies implemented to confront the health issues. This is the case for the aforementioned lockdowns that keep the individuals of a certain locality or sector confined and, therefore, reduce the economic activity. Therefore, it is possible to perceive a trade-off between the extension of lockdowns, and its impact on productivity, and the health outcomes observed; for instance, the longer lockdown period, the fewer intensive care units (ICU) required (or number of cases, or deaths). In terms of the objective for the decision makers, it is sought to minimize both the lockdowns duration and the health consequences \cite{Alvarez2020,Grigorieva2020}.

Several countries and cities have used  epidemiological indicators to activate, reinstate and release NPI such as lockdowns. 
In fact, in \cite{WHO-May2020}, the World Health Organization (WHO) has proposed  various public health  criteria to adjust public health and social measures in the context of COVID-19. The criteria are grouped into three domains: Epidemiology, Health system and Public health surveillance. Many of these criteria can be expressed in terms of the following indicators: effective reproduction number, active new cases, positivity rate, and hospitalization and intensive care unit (ICU) admissions due to COVID-19.

In this context, the government of San Francisco (US) has published on its website \cite{DataSF} various health indicators used to monitor the level of COVID-19  in the city and assess the ability of its health care system to respond to the pandemic. These health indicators are grouped into 5  areas: hospital system, cases, testing, contact tracing and personal protective equipment. The first three areas are characterized in terms of the following observations: COVID-19 hospitalizations, acute care beds available, ICU beds available, new cases per day per 100,000 residents, and tests collected per day. The last two areas consider the indicators related to contact tracing and personal protective equipment that are beyond the scope of our study. In other places, such as the city of Austin (US),  it was decided to track daily COVID-19 hospital admissions and daily total hospitalizations across the city and trigger the initiation and relaxation of lockdown periods when admissions cross predetermined thresholds \cite{Austin}. 

In the case of Chile, and in particular its capital, Santiago, which is our main case study, at the beginning of the outbreak, the government used indicators mainly based on the active cases (per 100,000 residents and per area) and on the available hospital beds to activate lockdowns \cite{minsalcuarentenas}. These lockdowns were dynamically applied to the most affected sectors of the cities. Since July 2020, a five-steps plan known as ``\emph{Paso a Paso}" (Step by Step), has been used to release and reinstate various NPI measures including lockdowns \cite{paso_a_paso}. Following the recommendations of the World Health Organization \cite{WHO-May2020} and of the  COVID-19 advisory council \cite{minuta2906}, an independent organism conceived to guide the Chilean Ministry of Public Health in the policies that will be implemented to face the COVID-19 outbreak; this plan now takes into account several dimensions of the pandemic through the monitoring of the following indicators: ICU beds available, effective reproductive number, new active cases, positivity rate, and some additional tracing and surveillance indicators. Thus, lockdowns are now release/reinstate for a given district of a city when the values of these indicators, computed for such districts, are below/above the predetermined thresholds. 

Given the aforementioned situation, our aim is to propose a framework based on control theory that permits to assess and compare NPI policies based on epidemiological indicators. 

It is important to note that the application of NPI strategies to mitigate the effects of the COVID-19 pandemic has been already modeled using optimal control techniques. For instance, in \cite{Djidjou-Demasse2020}, optimal control theory is used to explore the best strategy for implementation while waiting for the development of a vaccine. More specifically, they seek a solution minimizing deaths and costs due to the implementation of the control strategy itself. 

The same objective is studied in  \cite{Richard2020} by the means of an age-structured optimal control model. Specifically,  a model with a double continuous structure by host age and time since infection is proposed. In \cite{bonnans:hal-02558980}, the authors also propose an optimal control model with infection age. The main difference from the previous model is the  consideration of the peak value in the objective function that leads to a problem with state constraints. This work has generated an open access code, based on the optimal control toolbox BOCOP.   An SIR model is studied in  \cite{Angulo:2020} to derive a simple but mathematically rigorous criterion for designing optimal transitory NPIs. In particular, the authors found that reducing the reproduction number below one is sufficient but not necessary. This condition may be prescribed according to the maximum health services' capacity. 

In all of the cited articles, (optimal) control strategies are obtained numerically or analytically in closed-loop form (as in \cite{Angulo:2020}) requiring the assumption  that the state variable is perfectly observed. Moreover, these works consider that the effects of NPIs are perfectly controllable by the decision maker without taking into account that a decision regarding the application or release of an NPI may require time for the population to adapt. These considerations are important when the proposed trajectories are composed by bang-bangs, singular arcs or saturation of constraints, because in practice, these controls may be very difficult to implement. 

Additionally, since NPIs activation/release decisions based on the observation of indicators are made at specific given moments (for instance, lockdowns may be activated on Mondays for the whole week, or always at a given time of the day, etc.), they would not be faithfully represented through the aforementioned approaches. This is a key aspect when one of the objectives of the study is to provide recommendations to decision makers.

Thus, our first aim is to propose a framework for representing the decision-making  process related to the application of NPIs based on the observation of epidemiological indicators. The  proposed framework uses event-triggered controls in discrete time.  To the best of our knowledge, the notion of event-triggered control was first introduced in \cite{Arzen:1999} and allows to model update instants (in the decision, or control) by the violation of a condition depending on the state of the system, instead of assuming a constant or a periodic updating. As  mentioned above, the event-triggered decisions is precisely the framework such as the NPIs that are being applied in several places, because depending on whether or not some indicators are in violation of their respective thresholds, the NPIs are applied or lifted.

In recent years, event-triggered control techniques have attracted increasing interest with studies focusing on  the stabilization of event-triggered implementations and the existence of a minimum interevent time (see for instance \cite{8365671,Borgers:2014,Dolk:2017,Goebel:2009,Tabuada:2011,Tabuada:2007,Wakaiki:2020,ZhuLin:2019}). In our work, these topics are not addressed because it is assumed that the stabilization is not an issue, in the sense that we are not concerned with the final state of the system but rather we seek to be able to determine what type of indicator-based strategy is more suitable for use. Ensuring the minimum interevent time is not an issue in our framework, because we impose this minimum interevent time to avoid the application  and lifting (or vice versa) of an NPI in a short period of time, which is the practice we have observed in cited examples. Moreover, we consider the possibility  that decisions can be taken (apply or release an NPI) only at some prescribed times, for instance, only on Mondays, making the modeling of the decision-making process more realistic.

In addition to realistically representing the decision-making process based on the observation of indicators,  our approach allows the assessment and comparison of the use of indicators from a cost-benefit perspective. Indeed, based on the discussion given in this introduction section, we assume that a decision maker aims to obtain as few as possible of some outcomes that  are consequences of the NPIs. The peak of ICU beds demand and the total number of days in lockdown are two examples of such undesired outcomes.

Therefore, if the policy for activating  or releasing NPIs is triggered by some observation or, more specifically, when a given indicator (computed in terms of these observations) is above or below a certain threshold, we can link  all of these thresholds with their respective observations. 

Thus, considering the trigger threshold as a parameter, one can construct the curve of outcomes associated with a given indicator. This curve lies in the space of outcomes and it is called the \emph{trade-off curve}. This curve allows the identification of the trade-offs between different outcomes. Moreover, after computing these curves for several indicators, they can be compared to choose an  indicator that will determine  the policy. Indeed, for instance, if one of the curves dominates the others in the Pareto sense (that is, the curve is below the others), this immediately suggests that the corresponding indicator is the most suitable for the outcomes considered in the analysis. 

On the other hand, when no domination occurs between the curves, a decision can be taken by setting all but one outcomes as objectives  (viewed as upper bounds). This allows the determination of the best indicator as that with the lowest value for the outcome not considered as the objective. Thus, the proposed methodology  allows the comparison of the use of indicators to trigger either the application or release of NPIs.

The rest of this paper is organized as follows. In \cref{sec:preliminaries}, we introduce the proposed framework that we use to model the making decision process, that is, discrete-time control systems and event-triggered controls. In \cref{sec:tradeoff}, we define the trade-off curve associated with a trigger indicator and show how to compare these curves in order to choose the best indicator to activate/release NPIs. To illustrate  the proposed modeling framework and methodology for the assessment of the indicators, in \cref{sec:examples} we discuss two case-studies, corresponding to Chile and China, where different indicators are evaluated. Finally, some concluding remarks are stated  in \cref{sec:conclusions}.

\section{Mathematical framework based on event-triggered controls}\label{sec:preliminaries}

\subsection{Discrete-time control systems}

Given an initial time $\it \in \NN$, a horizon time $\horizon \in \NN$ (with $\horizon > \it$), an initial state $\ic \in \X$ and a finite sequence of controls $\u=(\control(t))_{t=\it}^{\horizon}=(\control(\it),\ldots,\control(\horizon))$, we consider the discrete-time control system:
\begin{equation}
\label{eq:system}
\state(t+1)=\dynamics(t,\state(t),\control(t)),\quad 
	t\in\titf{\it}{\horizon},\quad 
	\state(\it)=\state_0.	
  \end{equation}
  In this setting, $\dynamics:\titf{0}{N}\times\X\times\U\longrightarrow\X$ is the dynamics, $\X $ is a vector space, called the state space, and $\U$ is a compact set of a given vector space satisfying $0 \in \U$ (representing no action). The latter is called the control space. Here, when $p\leq q$,  $\titf{p}{q}$ denotes the collection of all integers between $p$ and $q$ (inclusive).
  
Several epidemiological models can be represented as system \eqref{eq:system}. For instance, the state $\state \in \X$ can model the number of individuals in different stages of a disease (e.g., SEIR compartmental models; See, for instance, \cite{castillo:2001} and references therein) in different populations (e.g., counties, cities, age-groups, etc.), and the control $\control \in \U$ can model the application of one or more NPIs to some groups of the population under study.
  
We denote by $\UU$ the collection of all possible controls; that is,
	$$\UU=\left\{\u=(\control(t))_{t=\it}^{\horizon}\middle|\ \control(\it),\ldots,\control(\horizon) \in \U \right\}\cong \U^{\horizon-\it+1}.$$

A solution of the control system \eqref{eq:system} associated with a control $\u\in\UU$ is an element of the space
$$\XX=\left\{\x=(\state(t))_{t=\it}^{\horizon+1}\middle|\ \state(\it),\ldots,\state(\horizon+1)\in\X\right\}\cong\X^{\horizon-\it+2},$$
that satisfies the initial time condition $\state(\it)=\ic$.

\subsection{Event-triggered controls}\label{sec:triggered}

For $\hist \in \NN$, with $0 \le \hist < \horizon -\it$, we  consider feedback controls based on the recent history of the state, that is, the state in the current time $t$ and in the previous $\hist$  times $t-1,\ldots,t-\hist$ (if $\hist > 0$). Therefore, for $t \in \titf{\it}{\horizon}$, we introduce the notation
$$\stateh_{\hist}(t)=(\state(t-\hist),\ldots,\state(t)) \in \X^{\hist+1}. $$
If $t-\hist < \it$, in the definition of $\stateh_{\hist}(t)$, we repeat the initial state $\ic$ in the first $\it -t +\hist $ components.

Motivated by the current practice, we consider $\stateh_{\hist}(t)$ instead of $\state(t)$ for defining our policies. Indeed, as mentioned in the introduction, the use of epidemiological indicators for triggering either the application, or the release, of NPIs considers the state of the disease in a window of time and not instantaneously.

To reduce the control variability,  we turn to the event-triggered
state feedback control based on the recent history of the state. For this purpose, we introduce an event-triggered set $\Set \subseteq  \X^{\hist+1}$ in order to determine the variations in the control law which in turn are produced by a transition of the recent history $\stateh_{\hist}(t)$, either from $\Set$ to its complement $\Set^c$, or vice versa.  To compute these transitions, we define the indicator (XOR) function $\feedbackf_{\Set}: \X^{\hist+1} \times \X^{\hist+1} \longrightarrow \{0,1\}$ by
$$\feedbackf_{\Set}(\stateh_{\hist},\tilde\stateh_{\hist})=\left\{\begin{array}{ll}
								0 & \mbox{ if } (\stateh_{\hist},\tilde\stateh_{\hist}) \in (\Set \times \Set) \cup (\Set^c \times \Set^c)\\[2mm]
								1  & \mbox{ if } (\stateh_{\hist},\tilde\stateh_{\hist}) \in (\Set \times \Set^c) \cup (\Set^c \times \Set) \end{array}\right.$$
for all $(\stateh_{\hist},\tilde\stateh_{\hist}) \in \X^{\hist+1} \times \X^{\hist+1}$. The event-triggered set $\Set$ will then  define the application ($\stateh_{\hist} \notin \Set$) or release ($\stateh_{\hist} \in \Set$) of an NPI.

Now, we define the following event-triggering mechanism
\begin{equation}\label{eq:controller0}
\control(t_0)=\left\{\begin{array}{ll}
							0& \mbox{ if } ~\state_{0} \in \Set\\[2mm]
							\control_{\rm ref} \in \U& \mbox{ if }~ \state_0\notin \Set
							\end{array}\right.
\end{equation}
and 
\begin{equation}\label{eq:controller}
\control(t)=\feedback(t-t_k,\stateh_{\hist}(t_k),\control(t_k)) \quad t_k < t \le t_{k+1},
\end{equation}
where the trigger times $t_k$  are updated according to rule \eqref{eq:times} given below, and  the controller $\feedback: \titf{0}{\horizon-\it} \times \X^{\hist+1}\times \U \longrightarrow \U$ is defined by 
\begin{equation}\label{eq:controllerk}
\feedback(\tau,\stateh_{\hist},\control)=\left\{\begin{array}{ll}
							\hat \control_-(\tau,\control) & \mbox{ if } ~\stateh_{\hist} \in \Set\\[2mm]
							\hat \control_+(\tau,\control)& \mbox{ if }~ \stateh_{\hist} \notin \Set ,
							\end{array}\right.
\end{equation}
with $\hat \control_-,~\hat \control_+: \titf{0}{\horizon-\it} \times \U \longrightarrow  \U$   reference time-varying controls  associated with the NPI to be applied. Functions $\hat \control_-(\cdot,\control)$ and $\hat \control_+(\cdot,\control)$  satisfy 
$\hat \control_-(0,\control)=\hat \control_+(0,\control)=\control$ for all $ \control \in \U$. The controller $\feedback(\tau,\stateh_{\hist},\control)$ is then defined (in \eqref{eq:controllerk}) through two types of policies depending on whether or not  $\stateh_{\hist}$ belongs to the event-triggered set  $\Set$. Thus, this models the application of a given NPI or the release of this measure. This policy, consisting of applying or not some NPIs,  is eventually time-varying on the intervals $\titf{t_k}{t_{k+1}}$, and can model a variable intensity or an adoption degree that changes with  time. For instance, in a population where an NPI has been applied for a long period of time, when the intervention is lifted,  some time is required to return to the activity levels (of contact rates, mobility, etc.) identical to those prior to the application of the NPI.  In addition, the controller takes into account the last control $\control$  applied in the previous interval $\titf{t_{k-1}}{t_{k}}$ with the aim of allowing (but not limiting to) smooth variations of the applied policies.

The trigger times $t_k$ in \eqref{eq:controller} are given by
\begin{equation}\label{eq:times}
t_{k+1}=t_k + \delay + \min\{t \ge 0~|~ \feedbackf_{\Set}(\stateh_{\hist}(t_k),\stateh_{\hist}(t_k+\delay+t))=1\} .
\end{equation}
In the above definition, we have considered a minimum fixed period of time $ \delay\in \titf{\hist}{\horizon - \it}$  during which the event-triggered control is applied, thus imposing a minimum  interevent time. The interevent time is assumed to be greater or equal to $\hist$ which is the time window where the state variable is considered in order to have enough time to observe the effects of the application of the NPI, or of its release.  On the other hand, the element
$$\stateh_{\hist}(t_k+\delay+t)=(\state(t_k+\delay+t-\hist),\ldots,\state(t_k+\delay+t))$$
in \eqref{eq:times} is the sequence of states obtained from
$$\state(\tau+1)=\dynamics(\tau,\state(\tau),\feedback(\tau-t_k,\stateh_{\hist}(t_k),\control(t_k))) \quad \tau \in \titf{t_k}{t_k+\delay +t-1}$$
with initial condition $\state(t_k)$.

More generally, trigger times in  \eqref{eq:times} can be defined as
$$t_{k+1}=t_k + \delay + \ddfunction\left(\min\{t \ge 0~|~ \feedbackf_{\Set}(\stateh_{\hist}(t_k),\stateh_{\hist}(t_k+\delay+t))=1\}\right)  ,$$
where the function $\ddfunction: \titf{0}{\horizon-\it} \longrightarrow \titf{0}{\horizon-\it}$ can be introduced in order to impose a change in the control policy in the prescribed periods. For instance, if the time unit is day and  the application (or release) of an NPI must start on a Monday, when $\it$ is a Monday, we take $\delay$ a multiple of seven and the function  $\ddfunction(t)=7 \ceil*{\frac{t}{7}}$, where $\ceil*{\alpha}$ is the smallest integer that is also greater or equal than $\alpha$.

To summarize, the event-triggered feedbacks considered in this paper are defined by:
\begin{enumerate}
\item \emph{Observation time window $\hist$}: period of time in which the state is observed for defining the closed-loop strategy;
\item \emph{Minimum time of implementation $\delay$}: period of time (with $\delay \ge \hist$) imposed for avoiding high variability in the control and in order to observe the effects of the applied measures;
 \item \emph{Event-triggered set $\Set$}: set of recent history of states for determining variations in the policy through transitions from $\Set$ to $\Set^c$ or in the inverse sense. The definition of set $\Set$ depends on $\hist$;

 \item \emph{Controller $\feedback$}: the function defining the control through mechanism described in \eqref{eq:controller0}, \eqref{eq:controller}, and \eqref{eq:controllerk}. The definition of  the controller $\feedback$ depends on functions  $\hat \control_-$ and $\hat \control_+$ in \eqref{eq:controllerk} and on $\hist$ and $\Set$.
  
\end{enumerate}
An  event-triggered feedback will be denoted by $\u=\trigger(\ic,\hist,\delay,\feedback,\Set) \in \UU$ in order to highlight that the associated control $\u$ is completely defined by the event-triggered mechanism (the four elements listed above) and the initial condition $\ic$ at time $\it$.

Hence, given the mathematical formalism described above, the main objective of this paper is to compare different even triggered controls assuming that  all of them have common controllers $\feedback$ with same functions $\hat \control_-$ and $\hat \control_+$ in \eqref{eq:controllerk} (that is, all of the policies apply the same NPI, and in the same manner). In addition,  since the initial state $\ic$ is fixed for all of our analyses, we simplify the notation of event-triggered feedbacks by $\u=\trigger(\hist,\delay,\Set)$.

\section{NPI strategies based on indicators and their trade-offs}\label{sec:tradeoff}

In this section, we specify the event-triggered set $\Set$ in the definition of event-triggered controls introduced in  \cref{sec:triggered}. The particular forms considered of the set $\Set$ are provided by an epidemiological indicator (observation). Thus, considering the space of the recent history of states $\X^{\hist+1}$ and an indicator $\indicator: \X^{\hist+1} \longrightarrow \R$, we will work with  event-triggered sets $\Set$ in the form
$$\Set=\Set(\hist,\indicator,\threshold)=\{\stateh_{\hist} \in \X^{\hist+1}~|~\indicator(\stateh_{\hist}) \le \threshold\} ,$$
where $\threshold \in \R$ is a threshold for the indicator $\indicator$.

If in the period of time $\titf{t-\hist}{t}$, we consider instantaneous observations   of the state, represented by the function $\observation: \X \longrightarrow \RR$, the following indicators $\indicator$ can be considered:

\begin{itemize}
\item[a)] the mean of instantaneous observations:
$$\indicator(\state_{\hist}(t))= \frac{1}{\hist +1} \sum_{\tau=t-\hist}^t \observation(\state(\tau)) ;$$
\item[b)] the mean of differences of instantaneous observations:
$$\indicator(\state_{\hist}(t))= \frac{1}{\hist +1} \sum_{\tau=t-\hist+1}^t (\observation(\state(\tau)) -\observation(\state(\tau-1)))= \frac{\observation(\state(t))-\observation(\state(t-\hist))}{\hist +1};$$
\item[c)] variation rate of instantaneous observations:
$$\indicator(\state_{\hist}(t))= \frac{ \observation(\state(t))-\observation(\state(t-\hist)) }{\observation(\state(t-\hist)) };$$

\item[d)] variation rate of differences of instantaneous observations:
$$\indicator(\state_{\hist}(t))=  \frac{1}{\hist +1} \sum_{\tau=t-\hist}^t  \frac{ \observation(\state(t))-\observation(\state(t-1)) }{\observation(\state(t-1)) }.$$
\end{itemize}
The specific instances of observations $\observation$ can be the number of  infected people \cite{minuta2906,minsalcuarentenas,paso_a_paso,DataSF,WHO-May2020}, number of deaths associated with COVID-19 \cite{WHO-May2020}, excess mortality due to pneumonia \cite{WHO-May2020}, number of tests \cite{DataSF}, positivity rate \cite{paso_a_paso,WHO-May2020}, COVID-19 total hospitalizations \cite{Austin,paso_a_paso,DataSF,WHO-May2020} and the number of people in ICU beds \cite{paso_a_paso,DataSF,WHO-May2020}, among others. Additionally, the formalism introduced above permits to consider more elaborated indicators $\indicator$ computed from $\state_{\hist}(t)$, such as the effective reproductive number \cite{paso_a_paso,WHO-May2020}. 

We update the notation of an event-triggered feedback  by $\u=\trigger(\hist,\delay,\indicator,\threshold)$, noticing that the event-triggered set $\Set$ now depends on indicator $\indicator$ and threshold $\threshold$.

\subsection{Trade-offs of event-triggered policies based on indicators}\label{sec:outcomes}

For a fixed  initial state $\ic \in \X$ and a given control sequence $\u \in \UU$,  we assume that $m \in \NN$ ($m \ge 2$) outcome  measures are observed. These outcomes are represented by a function $\output: \X \times \UU \longrightarrow \R^m$ and, in our convention, we want these to be as low as possible. Examples of these outcomes are the peak demand of ICU beds, the total number of deaths (directly attributed to the COVID-19 outbreak) and the total number of days of lockdown (i.e., shelter-in-place), or of any other NPI. The first two outcomes are related to the health cost of to the pandemic, while the third-one is a proxy for the economic and societal cost of the policy. Intuitively, there should be a trade-off between the different outcomes, in the sense that a lower outcome of one kind (e.g., peak demand of ICU beds) will lead to a greater outcome of the other kind  (e.g., total time of NPI), as it is proven in \cite{Angulo:2020}. Of course, one can take a group of outcomes where the mentioned trade-off does not exist (e.g., peak demand of ICU beds and the  total number of deaths), but even in these cases it is interesting to quantify the links between these outcomes. Therefore, we propose to compute the outcomes given by an event-triggered feedback  $\u=\trigger(\hist,\delay,\indicator,\threshold)$ when the threshold $\threshold$ varies. With this information, decision makers can observe the links between the outcomes provided by the policy based on the indicator $\indicator$ and, for instance, to evaluate the marginal variation of one outcome when other outcome decreases/increases.

For a given event-triggered feedback  $\u=\trigger(\hist,\delay,\indicator,\threshold)$, we define \emph{the trade-off curve} associated with indicator $\indicator$, and parametrized by thresholds $\threshold$, as follows
\begin{equation}\label{eq:curve}
\tradeoff(\ic,\hist,\delay,\indicator)=\{\output(\ic,\u) ~|~\u=\trigger(\hist,\delay,\indicator,\threshold);~~\threshold \in \thresholdd_{\indicator}\} \subseteq \R^m ,
\end{equation}
where $\thresholdd_{\indicator} \subseteq \R$ is the domain for the thresholds associated with indicator $\indicator$.

\begin{figure}[ht]
\begin{center}
\psset{xunit=.5cm, yunit=.5cm, linewidth=.2pt}
\begin{pspicture}(1,-2)(12,13)
\psaxes[linecolor=black,  arrows=->, linewidth=1.5pt,Ox=,Oy=,Dx=20,Dy=2,dx=2.5,dy=3,labels=none,tickstyle=bottom,ticksize=.2](0,0)(15,12)
\pscurve[linewidth=1pt, linecolor=red,fillstyle=none,fillcolor=palegreen](.3,11.5)(1,6)(3,3)(14,1)
\psline[linewidth=1pt, linecolor=white,fillstyle=none,fillcolor=palegreen](.28,11.5)(13.95,11.5)(13.95,1)

\rput(-0.9,12.1) {$\outputc_2$}
\rput(15.5,-0.7) {$\outputc_1$}
\rput(5,4) {\textcolor{red}{$\tradeoff(\ic,\hist,\delay,\indicator)$}}
\psline[linewidth=1pt, linecolor=black,linestyle=dashed](6.3,1.7)(6.3,0)
\psline[linewidth=1pt, linecolor=black,linestyle=dashed](6.3,1.7)(0,1.7)
\rput(6.3,1.7){\textcolor{red}{$\bullet $}}
\rput(6.3,-.6) {$\outputc_1(\threshold_1)$}
\rput(-1.4,1.7) {$\outputc_2(\threshold_1)$}
\psline[linewidth=1pt, linecolor=black,linestyle=dashed](2,4)(2,0)
\psline[linewidth=1pt, linecolor=black,linestyle=dashed](2,4)(0,4)
\rput(2,4){\textcolor{red}{$\bullet $}}
\rput(2,-.6) {$\outputc_1(\threshold_2)$}
\rput(-1.4,4) {$\outputc_2(\threshold_2)$}

\end{pspicture}
 \end{center}
 \caption{Illustration of a trade-off curve $\tradeoff(\ic,\hist,\delay,\indicator)$ parametrized by thresholds $\threshold$  when the number of outputs is $m=2$. }\label{fig:tradeoff}
\end{figure}
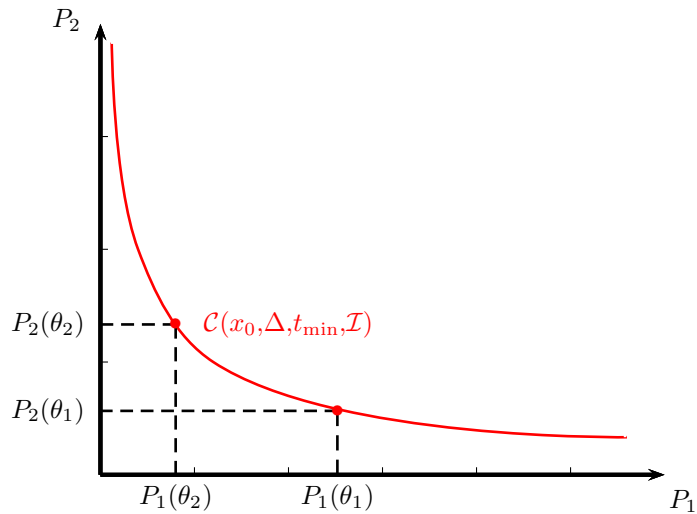

In \cref{fig:tradeoff}, we illustrate the  trade-off curve $\tradeoff(\ic,\hist,\delay,\indicator)$ when the observed outcomes are $m=2$, that is, $\output(\ic,\u)=(\outputc_1(\ic,\u),\outputc_2(\ic,\u))$ with $\u=\trigger(\hist,\delay,\indicator,\threshold)$. In this figure,  we use the simplified notation $ \output(\ic,\u)=(\outputc_1(\theta),\outputc_2(\theta)) \in \R^2$, for highlighting the dependence of the outcomes on thresholds. In this curve, we represent two outcomes exhibiting a trade-off (i.e., greater value of one outcome leads to a lower values of the other outcome).
 
\subsection{Comparing  event-triggered policies based on indicators}\label{sec:comparison}

To compare the efficiency of  two  event-triggered feedbacks based on  two different epidemiological indicators $\indicator_a$ and $\indicator_b$, we propose to compute the associated trade-off curves described in the previous section. If the policies based on indicators $\indicator_a$ and $\indicator_b$ consider observation time windows $\hist_a$ and $\hist_b$, respectively, and minimum interevent times  $\delay^a$ and  $\delay^b$, respectively, then we compare the trade-off curves $\tradeoff_a=\tradeoff(\ic,\hist_a,\delay^a,\indicator_a)$ and $\tradeoff_b=\tradeoff(\ic,\hist_b,\delay^b,\indicator_b)$ defined in \eqref{eq:curve}.

Recall that the outcomes are represented by the function $\output: \X \times \UU \longrightarrow \R^m$, with $m \ge 2$, and  the aim is to have the $m$ outcomes $\output(\ic,\u)=(\outputc_1(\ic,\u),\ldots,\outputc_m(\ic,\u))$ to be as low as possible. Then, the calculation of the trade-off curves $\tradeoff_a$ and $\tradeoff_b$, consists of the evaluation of function $\output$ for  event-triggered controls $\u_a(\threshold_a)=\trigger(\hist_a,\delay^a,\indicator_a,\threshold_a)$ and $\u_b(\threshold_b)=\trigger(\hist_b,\delay^b,\indicator_b,\threshold_b)$ with different values of thresholds $\threshold_a$ and $\threshold_b$.

\begin{figure}[ht]
\begin{center}
\psset{xunit=.5cm, yunit=.5cm, linewidth=.2pt}

\begin{pspicture}(1,-2)(12,13)
\psaxes[linecolor=black,  arrows=->, linewidth=1.5pt,Ox=,Oy=,Dx=20,Dy=2,dx=2.5,dy=3,labels=none,tickstyle=bottom,ticksize=.2](0,0)(15,12)

\pscurve[linewidth=1pt, linecolor=blue,fillstyle=none](.3,11.5)(2,6)(3,5)(14,2)
\psline[linewidth=1pt, linecolor=white,fillstyle=none](.28,11.5)(13.95,11.5)(13.95,1)

\pscurve[linewidth=1pt, linecolor=red,fillstyle=none](1,11.5)(3,6)(6,3)(14,1)
\psline[linewidth=1pt, linecolor=white,fillstyle=none](1.5,11.5)(13.95,11.5)(13.95,3)
\rput(5.5,10) {\textcolor{red}{$\tradeoff_a=\tradeoff(\ic,\hist_a,\delay^a,\indicator_a)$}}
\rput(12,4) {\textcolor{blue}{$\tradeoff_b=\tradeoff(\ic,\hist_b,\delay^b,\indicator_b)$}}

\rput(-0.9,12.1) {$\outputc_2$}
\rput(15.5,-0.7) {$\outputc_1$}

\psline[linewidth=1pt, linecolor=black,linestyle=dashed](2,8)(2,0)
\psline[linewidth=1pt, linecolor=black,linestyle=dashed](2,8)(0,8)
\psline[linewidth=1pt, linecolor=black,linestyle=dashed](2,6)(0,6)

\psline[linewidth=1pt, linecolor=black,linestyle=dashed](11,2.5)(11,0)
\psline[linewidth=1pt, linecolor=black,linestyle=dashed](11,2.5)(0,2.5)
\psline[linewidth=1pt, linecolor=black,linestyle=dashed](11,1.5)(0,1.5)

\rput(2,8.05){\textcolor{red}{$\bullet $}}
\rput(2,6){\textcolor{red}{$\bullet $}}

\rput(11,2.5){\textcolor{red}{$\bullet $}}
\rput(11,1.5){\textcolor{red}{$\bullet $}}
\rput(11,-.6) {$\hat \outputc_1$}
\rput(2,-.6) {$\hat \outputc_1$}
\rput(-1,2.5) {$\outputc_2^b$}
\rput(-1,1.5) {$\outputc_2^a$}

\rput(-1,6) {$\outputc_2^b$}
\rput(-1,8) {$\outputc_2^a$}
\end{pspicture}
 \end{center}
 \caption{Illustration of the trade-off curves $\tradeoff_a$ and $\tradeoff_b$ corresponding to event-triggered feedbacks based on indicators $\indicator_a$ and $\indicator_b$, when the number of outcomes $m=2$ and there is no domination.}\label{fig:domina2}
\end{figure}
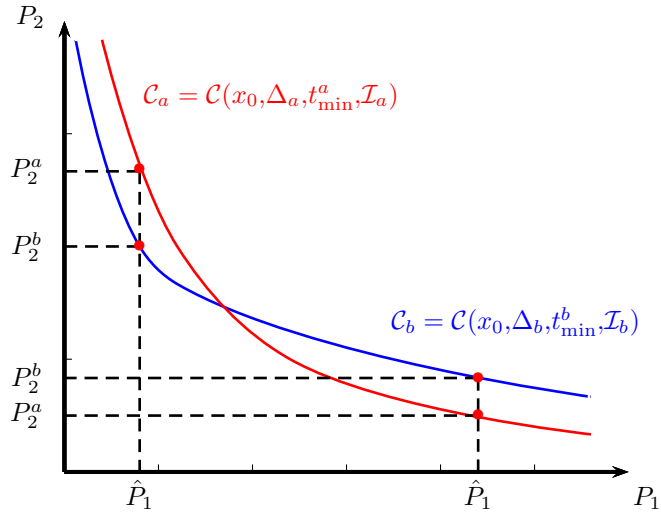

If a decision maker considers as objective (maximal bounds) the first $m-1$ outcomes $\outputcv{1}{m-1}=(\hat \outputc_1,\ldots, \hat \outputc_{m-1})$ (e.g., peak demand of ICU beds), the corresponding $m$-th outcomes $\outputc_m^a$ and $\outputc_m^b$ (e.g., total time of NPI) obtained using event-triggered feedbacks based on indicators $\indicator_a$ and $\indicator_b$ can be computed such that
$$(\outputcv{1}{m-1},\outputc_m^a) \in \tradeoff_a ~~~\mbox{ and } ~~~(\outputcv{1}{m-1},\outputc_m^b) \in \tradeoff_b .$$
Thus,  the decision maker should consider the policy based on $\indicator_a$ or $\indicator_b$ comparing $\outputc_m^a$ and $\outputc_m^b$ decision depending on the given objective outcomes $\outputcv{1}{m-1}$. In \cref{fig:domina2}, we show two trade-off curves for $m=2$. In this figure, for small objective values of outcome $\hat \outputc_1$, it is better to use a policy based on indicator $\indicator_b$. Otherwise, for larger values of $\hat \outputc_1$, it is recommended to use a policy based on $\indicator_a$.

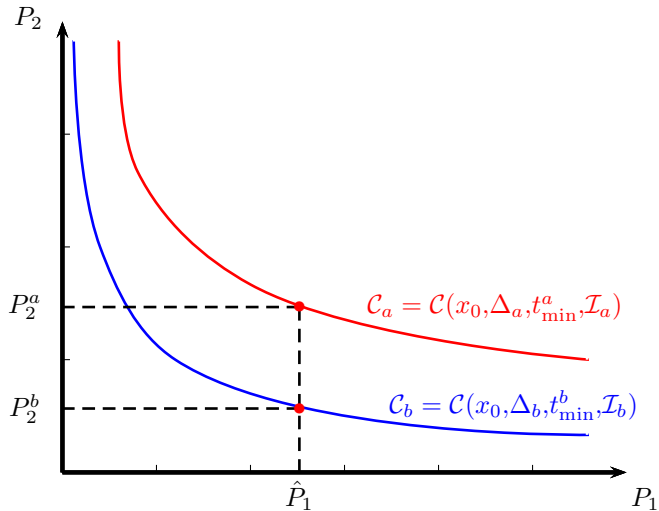
\begin{figure}[H]
\begin{center}
\psset{xunit=.5cm, yunit=.5cm, linewidth=.2pt}

\begin{pspicture}(1,-2)(12,13)
\psaxes[linecolor=black,  arrows=->, linewidth=1.5pt,Ox=,Oy=,Dx=20,Dy=2,dx=2.5,dy=3,labels=none,tickstyle=bottom,ticksize=.2](0,0)(15,12)

\pscurve[linewidth=1pt, linecolor=blue,fillstyle=none](.3,11.5)(1,6)(3,3)(14,1)
\psline[linewidth=1pt, linecolor=white,fillstyle=none](.28,11.5)(13.95,11.5)(13.95,1)

\pscurve[linewidth=1pt, linecolor=red,fillstyle=none](1.5,11.5)(2,8)(5,5)(14,3)
\psline[linewidth=1pt, linecolor=white,fillstyle=none](1.5,11.5)(13.95,11.5)(13.95,3)
\rput(11.5,4.4) {\textcolor{red}{$\tradeoff_a=\tradeoff(\ic,\hist_a,\delay^a,\indicator_a)$}}
\rput(12,1.8) {\textcolor{blue}{$\tradeoff_b=\tradeoff(\ic,\hist_b,\delay^b,\indicator_b)$}}

\rput(-0.9,12.1) {$\outputc_2$}
\rput(15.5,-0.7) {$\outputc_1$}

\psline[linewidth=1pt, linecolor=black,linestyle=dashed](6.3,4.4)(6.3,0)
\psline[linewidth=1pt, linecolor=black,linestyle=dashed](6.3,4.4)(0,4.4)
\psline[linewidth=1pt, linecolor=black,linestyle=dashed](6.3,1.7)(0,1.7)

\rput(6.3,1.7){\textcolor{red}{$\bullet $}}
\rput(6.3,4.4){\textcolor{red}{$\bullet $}}
\rput(6.3,-.6) {$\hat \outputc_1$}
\rput(-1,1.7) {$\outputc_2^b$}
\rput(-1,4.4) {$\outputc_2^a$}

\end{pspicture}
 \end{center}
 \caption{Illustration of trade-off curves $\tradeoff_a$ and $\tradeoff_b$ corresponding to event-triggered feedbacks based on indicators $\indicator_a$ and $\indicator_b$, when the number of outcomes is $m=2$ and one curve dominates the other.}\label{fig:domina}
\end{figure}

 If $\outputc_m^b < \outputc_m^a$ for all objective outcomes $\outputcv{1}{m-1}$, it means that the curve $\tradeoff_b$ dominates curve $\tradeoff_a$ and therefore, the policy based on $\indicator_b$ is always better. This situation is illustrated in \cref{fig:domina}  for $m=2$.

Finding a control $\u \in \UU$ that minimizes  the $m$ outcomes given by the vector $\output(\ic,\u)=(\outputc_1(\ic,\u),\ldots,\outputc_m(\ic,\u))$ is evidently a multicriteria optimization problem. Therefore, we only propose to compare single outcomes obtained from two policies  when the other ($m-1$) outcomes are fixed and treated as objectives. This avoids the need to consider outcomes of different nature and expressed in different units in a single objective function.  On the other hand, to compare two policies when several outcomes are simultaneously considered (i.e., not fixed) is technically possibly, although it is not straightforward.

\section{Examples}\label{sec:examples}

\subsection{Case study: Metropolitan Region, Chile}\label{sec:chile}

For this example, we  use a  discrete-time compartmental model  where the population of Metropolitan Region, Chile (Santiago)  $N= 7,112,808$ (according to the Chilean Institute of Statistics INE\footnote{\url{https://www.ine.cl/estadisticas/sociales/censos-de-poblacion-y-vivienda/poblacion-y-vivienda}}), assumed to be constant and isolated, is distributed into eight groups corresponding to different disease stages. Susceptible (denoted by $\susceptible$) are individuals not infected by the disease but who can be infected by the virus. Exposed (denoted by $\exposed$) are those in the incubation time after being infected. In this stage, they do not have symptoms but can infect other people with a lower probability than those in the infectious group described below. Mild infected or subclinical (denoted by $\infectedU$) is an infected population that can also infect other people. In this stage, they are asymptomatic or show mild symptoms. They are not detected and thus are not reported by authorities. At the end of this stage, they pass directly to the recovered state. Infected (denoted by $\infected$) are infected citizens that can infect other people. They develop symptoms and are detected and reported by authorities. They can recover or enter some hospitalized state. Recovered (denoted by $\recovered$) is a population that survives the illness, is no longer infectious, and has developed immunity to the disease. Hospitalized (denoted by $\hospitalized$) are patients hospitalized in basic facilities.  After this stage, hospitalized can recover or get worse and use an ICU bed or die. Hospitalized in ICU beds (denoted by $\hospitalizedC$) are patients in critical care. 
People in these two last stages may infect other people; however, in our posterior analysis, we will neglect this source of infection because we assume that hospitalized individuals are highly isolated. After leaving the ICU bed, they either come back to the basic facilities in the hospital or die. The term dead (denoted by $\dead$) accounts for people who did not survive the disease.

The evolution of the previous state variables is described by the following discrete-time system, with states $\state=(\susceptible,\exposed,\infectedU,\infected,\recovered,\hospitalized,\hospitalizedC,\dead) \in \X = \R^8_+$:

\begin{equation}
\label{eq:model}
\left\{\begin{array}{lcl}

\susceptible(t+1) &=& \susceptible(t) - \rateg(\state(t))\susceptible(t) 
 \\[3mm]

\exposed(t+1) &=& \exposed(t)+\rateg(\state(t)) \susceptible(t)- \duration_{\exposed}\exposed(t) \\[3mm]

\infectedU(t+1) &=& \infectedU(t)+(1-\fraction_{\exposed \infected})\duration_{\exposed} \exposed(t) -\duration_{\infectedU}\infectedU(t)\\[3mm]

\infected(t+1) &=& \infected(t)+\fraction_{\exposed \infected}\duration_{\exposed} \exposed(t) -\duration_{\infected}\infected(t)\\[3mm]

\recovered(t+1) &=& \recovered(t)+\duration_{\infectedU}\infectedU(t)+ \fraction_{\infected\recovered} \duration_{\infected} \infected(t) + \fraction_{\hospitalized\recovered}\duration_{\hospitalized} \hospitalized(t)\\[3mm]

\hospitalized(t+1) &=& \hospitalized(t)+(1- \fraction_{\infected\recovered})\duration_{\infected} \infected(t) +(1-\fraction_{\hospitalizedC\dead}) \duration_{\hospitalizedC} \hospitalizedC(t)- \duration_{\hospitalized}\hospitalized(t)\\[3mm]

\hospitalizedC(t+1) &=& \hospitalizedC(t)+(1-\fraction_{\hospitalized\recovered}-\fraction_{\hospitalized\dead}) \duration_{\hospitalized} \hospitalized(t) - \duration_{\hospitalizedC} \hospitalizedC(t)\\[3mm]

\dead(t+1) &=& \dead(t)+\fraction_{\hospitalized \dead} \duration_{\hospitalized} \hospitalized(t)+ \fraction_{\hospitalizedC \dead} \duration_{\hospitalizedC} \hospitalizedC(t) ,

\end{array}\right.
\end{equation}
where natural births and deaths are not considered because their effects are negligible.

The above system is a type of discrete-time SEIR (or SEIRHD) model. It aims to better describe an outbreak where part of the population has been infected by a virus, with many people with no symptoms or only mild symptoms. It was found that this is the case for SARS-CoV-2, as has been reported in several papers in the literature \cite{imperial,ivorra:2020,magal}. The model we use can be extended to consider age ranges and evaluate other indicators or policies such as the schools reopening or the progressive lifting of quarantines according to age groups as explored in \cite{castillo-levin,Zhao:2020}. However, to illustrate the proposed methodology, we have only considered one age class. 

\begin{figure}[ht]
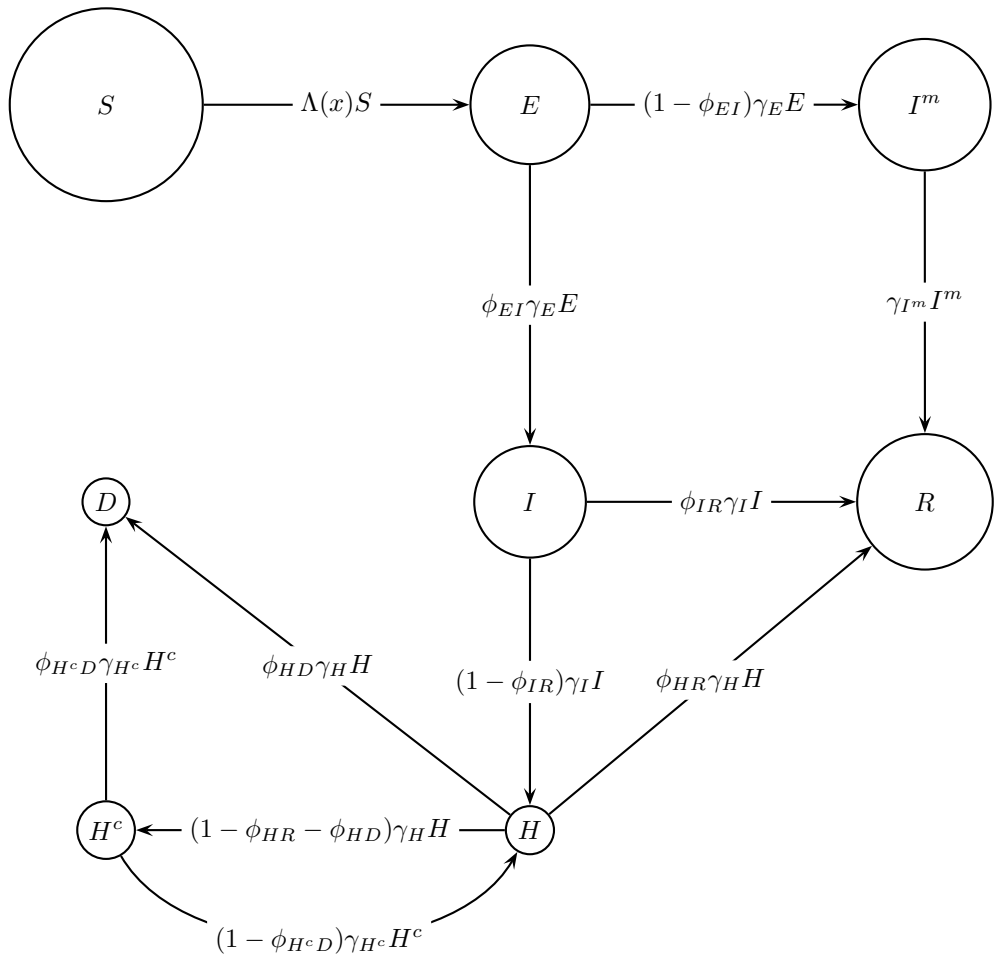

\begin{equation*}
\psmatrix[colsep=3.5cm,rowsep=3cm,mnode=circle,arrowscale=1.5,radius=15pt]
\textcolor{white}{Xx}\textcolor{white}{Xx}\susceptible \textcolor{white}{Xx}\textcolor{white}{xX}&\textcolor{white}{Xx}\exposed\textcolor{white}{Xx}&\textcolor{white}{Xx}\infectedU\textcolor{white}{Xx}\\
\dead&\textcolor{white}{Xx}\infected\textcolor{white}{Xx}&\textcolor{white}{X} \textcolor{white}{X} \recovered  \textcolor{white}{X}  \textcolor{white}{X}\\
\hospitalizedC & \hospitalized& [mnode=none]
\ncline{->}{1,1}{1,2}
\ncput*{\rateg(\state)\susceptible  }
\ncline{->}{1,2}{1,3}
\ncput*{(1-\fraction_{\exposed \infected})\duration_{\exposed} \exposed}
\ncline{->}{1,2}{2,2}
\ncput*{\fraction_{\exposed \infected}\duration_{\exposed} \exposed}
\ncline{->}{2,2}{2,3}
\ncput*{\fraction_{\infected \recovered}\duration_{\infected} 	\infected}
\ncline{->}{1,3}{2,3}
\ncput*{\duration_{\infectedU}\infectedU }
\ncline{->}{2,2}{3,2}
\ncput*{(1-\fraction_{\infected \recovered})\duration_{\infected}\infected}
\ncline{->}{3,2}{3,1}
\ncput*{(1-\fraction_{\hospitalized \recovered}-\fraction_{\hospitalized\dead})\duration_{\hospitalized}\hospitalized}
\ncline{->}{3,2}{2,3}
\ncput*{\fraction_{\hospitalized \recovered} \duration_{\hospitalized}\hospitalized }
\ncarc[arcangle=-60]{->}{3,1}{3,2}
\ncput*{(1-\fraction_{\hospitalizedC\dead})\duration_{\hospitalizedC}\hospitalizedC}
\ncline{->}{3,1}{2,1}
\ncput*{\fraction_{\hospitalizedC\dead}\duration_{\hospitalizedC}\hospitalizedC}
\ncline[arcangle=120]{->}{3,2}{2,1}
\ncput*{\fraction_{\hospitalized\dead}\duration_{\hospitalized}\hospitalized}
\endpsmatrix\end{equation*}
\vspace{1cm}
\caption{{\bf Structure of the mathematical model for the COVID-19 dynamics in an isolated city (Metropolitan Region, Chile). }
Each circle represents a specific group. Susceptible individuals ($\susceptible$), and different disease states:  exposed ($\exposed$), mild  infected ($\infectedU$), infected ($\infected$), recovered ($\recovered$), hospitalized ($\hospitalized$), hospitalized in ICU beds ($\hospitalizedC$), and dead ($\dead$). }
\label{fig:model}
\end{figure}

The structure of this mathematical model, along with the transitions among the different stages, is shown in  \cref{fig:model}. 

We proceed now to describe the main components of system \eqref{eq:model}.  First, we consider the time interval $\titf{\it}{\horizon}$, with days as the unit of time,  where $t_0$ is the initial time from which we start the assessment  of different policies (September 21, 2020) and $\horizon$ is the  time horizon. In this example, the time horizon is the simulation time and must be large enough for all of the considered policies (NPI triggered by indicators) to be able to  reduce the number of infected people to zero (September 21, 2025), without considering other interventions. 
 
The contagion rate (see the first and second equation in \eqref{eq:model}) is
$$
\rateg(\state(t))=\frac{1}{(N-\dead(t))}
(\rate_{\exposed} \exposed(t)  +\rate_{\infectedU} \infectedU(t)+\rate_{\infected} \infected(t) + \rate_{\hospitalized}\hospitalized(t)+\rate_{\hospitalizedC}\hospitalizedC(t)) ,
$$
where the specific rates $\rate_{\exposed}, \rate_{\infectedU}, \rate_{\infected}, \rate_{\hospitalized},  \rate_{\hospitalizedC}$ can be written as the product of the contagion probabilities ($\prob_{\exposed}$, $\prob_{\infectedU}$, $\prob_{\infected}$, $\prob_{\hospitalized}$, and $\prob_{\hospitalizedC}$) and the reference values of contact rates ($\crref_{\exposed}$, $\crref_{\infectedU}$, $\crref_{\infected}$, $\crref_{\hospitalized}$, and $\crref_{\hospitalizedC}$), each value associated with a contagious stage of the disease (indicated in the subindex). Therefore, one has $\rate_{\exposed}=\prob_{\exposed}\crref_{\exposed}$, ~$ \rate_{\infectedU}=\prob_{\infectedU}\crref_{\infectedU}$, ~ $\rate_{\infected}=\prob_{\infected}\crref_{\infected}$, ~
$\rate_{\hospitalized}=\prob_{\hospitalized}\crref_{\hospitalized}$,~ and 
$\rate_{\hospitalizedC}=\prob_{\hospitalizedC}\crref_{\hospitalizedC}$.

Our goal is to influence the contagion rate $\rateg(\state(t))$ with a control policy $\control(t)$, representing the application of NPI, in this example, a lockdown, having as effect a reduction in the contact rates.  We denote by $\hat \rateg(\state(t),\control(t))$ the controlled contagion rate and  we define controlled specific contagion rates $\hat{\rate}_{\exposed}, \hat{\rate}_{\infectedU}, \hat{\rate}_{\infected}, \hat{\rate}_{\hospitalized}, \hat{\rate}_{\hospitalizedC}$ as follows
\begin{equation}
\label{eq:controlledrates}
\left\{\begin{array}{lclcl}

\hat{\rate}_{\exposed}(\control_{\exposed}) &=& (1-\control_{\exposed})\rate_{\exposed} &=&   (1-\control_{\exposed})\prob_{\exposed}\crref_{\exposed} \\[3mm]

\hat{\rate}_{\infectedU}(\control_{\infectedU}) &=& (1-\control_{\infectedU})\rate_{\infectedU} &=& (1-\control_{\infectedU})\prob_{\infectedU}\crref_{\infectedU}\\[3mm]

\hat{\rate}_{\infected}(\control_{\infected}) &=& (1-\control_{\infected})\rate_{\infected} &=& (1-\control_{\infected})\prob_{\infected}\crref_{\infected}\\[3mm]

\hat{\rate}_{\hospitalized}(\control_{\hospitalized}) &=& (1-\control_{\hospitalized})\rate_{\hospitalized} &=&  (1-\control_{\hospitalized})\prob_{\hospitalized}\crref_{\hospitalized} \\[3mm]

\hat{\rate}_{\hospitalizedC}(\control_{\hospitalizedC}) &=& (1-\control_{\hospitalizedC})\rate_{\hospitalizedC} &=&  (1-\control_{\hospitalizedC})\prob_{\hospitalizedC}\crref_{\hospitalizedC}

\end{array}\right.
\end{equation}
where $u=(\control_{\exposed},\control_{\infectedU},\control_{\infected},\control_{\hospitalized},\control_{\hospitalizedC}) \in [0,1]^5$ is the initial vector of control variables that later will be reduced.

We shall assume that $\crref_{\hospitalized} = \crref_{\hospitalizedC}=0$ because we suppose  that the hospitalized patients are highly isolated. Therefore,  $ \control_{\hospitalized}$ and $ \control_{\hospitalizedC}$ are no longer considered as control variables. This approach is also used in \cite{singapure}. Additionally, we assume that the policy for the exposed and mild-infected individuals is the same and thus $\control_{\exposed} = \control_{\infectedU}$, which will be denoted $\control $ for simplicity. This means that  NPI have the same effects (in terms of the contact rates reduction) for the people in the incubation stage and for the infected population with mild symptoms because these populations are not detected, so that NPI such as  lockdowns are assumed to affect them similarly. Finally, we also assume that NPI do not have additional effects in terms of the contact rates reduction, in the infected population with symptoms (and then detected and reported by authorities) because we assume that these individuals are already isolated. Thus, we set $\control_{\infected}=1-\hat \delta$  with $\hat \delta \in (0,1)$ representing the fraction of nonisolated infected (detected) population. Our vector of control variables is then reduced to a single input $\control \in [0,1]$ and the (controlled) contagion rate becomes
\begin{equation}\label{eq:controlledrate}
\hat \rateg(\state(t),\control(t))=\frac{1}{(N-\dead(t))}\left(
(1-\control(t))
(\rate_{\exposed} \exposed(t)  +\rate_{\infectedU} \infectedU(t))+\hat{\delta}\rate_{\infected} \infected(t)\right) .
\end{equation}
Furthermore, we set $1-\hat \delta$ as the upper bound for the control variable, thus representing that  NPI cannot isolate people to a greater extent than the isolation already imposed on the infected and detected population. Therefore, the single control variable $\control$ belongs to the control space $\U=[0,1-\hat \delta]$. In the following, we set $\hat \delta=0.2$.

Now we can write the (autonomous) dynamics $\dynamics:\X\times\U\longrightarrow\X$, where the state space is $\X=\R^8_+$, defining the discrete-time control systems of this example as in \eqref{eq:system} as follows

\begin{equation}
\label{eq:modelf}
\dynamics(\state,\control)=\dynamics(\susceptible,\exposed,\infectedU,\infected,\recovered,\hospitalized,\hospitalizedC,\dead,\control) =
\left(\begin{array}{c}

\susceptible -\hat \rateg(\state,\control)\susceptible 
 \\[3mm]

 \exposed+ \hat \rateg(\state,\control)\susceptible - \duration_{\exposed}\exposed \\[3mm]

 \infectedU+(1-\fraction_{\exposed \infected})\duration_{\exposed} \exposed -\duration_{\infectedU}\infectedU\\[3mm]

\infected+\fraction_{\exposed \infected}\duration_{\exposed} \exposed -\duration_{\infected}\infected\\[3mm]

 \recovered+\duration_{\infectedU}\infectedU+ \fraction_{\infected\recovered} \duration_{\infected} \infected + \fraction_{\hospitalized\recovered}\duration_{\hospitalized} \hospitalized\\[3mm]

 \hospitalized+(1- \fraction_{\infected\recovered})\duration_{\infected} \infected +(1-\fraction_{\hospitalizedC\dead}) \duration_{\hospitalizedC} \hospitalizedC- \duration_{\hospitalized}\hospitalized\\[3mm]

 \hospitalizedC+(1-\fraction_{\hospitalized\recovered}-\fraction_{\hospitalized\dead}) \duration_{\hospitalized} \hospitalized - \duration_{\hospitalizedC} \hospitalizedC\\[3mm]

 \dead+\fraction_{\hospitalized \dead} \duration_{\hospitalized} \hospitalized+ \fraction_{\hospitalizedC \dead} \duration_{\hospitalizedC} \hospitalizedC ,

\end{array}\right) ,
\end{equation}
where the controlled contagion rate $\hat \rateg(\state,\control)$  is given by \eqref{eq:controlledrate}.  The rest of parameters defining the above dynamics are described in \cref{app:chile} where we describe the calibration procedure for this example.

Now we describe the elements for defining the event-triggered controls introduced in \cref{sec:triggered}. The considered \emph{observation time window} and the \emph{minimum time of implementation} (or minimum interevent time) will be $\hist=\delay=14$ days. The event-triggering mechanism, as described in \eqref{eq:controller0} and  \eqref{eq:controllerk}, is defined by $\control_{\rm ref}=1-\hat \delta$ (associated with $\control(\it)$ in \eqref{eq:controller0})  and  as \emph{controllers} $\hat \control_-,~\hat \control_+: \titf{0}{\horizon-\it} \times \U \longrightarrow  \U$ (see  \eqref{eq:controllerk}),
we consider max-linear functions  of the form
\begin{eqnarray*}
\hat \control_-(\tau,\control)&=&\max\left\{0, \control\left(1 - \frac{\tau}{\delay}\right)\right\}\\[2mm]
\hat \control_+(\tau,\control)&=&\min\left\{1 -\hat \delta, \control\left(1 - \frac{\tau}{\delay}\right)  + \left(\frac{1-\hat \delta}{\delay}\right) \tau\right\} .
\end{eqnarray*}
These controllers permit to monotonically decrease/increase the strength of the NPI in  $\delay$ days until saturating it. We have decided to consider that our control saturates in $\delay$ days, coinciding with the minimum interevent time assuming that at the end of  this period of time, the effect of NPI (or their release) is fully accomplished, producing the maximal isolation ($\control = 1-\hat \delta$) or reaching the normal contact rates equal to those before  the pandemic ($\control = 0$). 

To define trigger indicators $\indicator: \X^{\hist+1} \longrightarrow \R$ and then the event-triggered sets $\Set$, as in \cref{sec:tradeoff}, we consider the following instantaneous observations of the state  $\observation: \X \longrightarrow \RR$:
\begin{enumerate}[label=(\subscript{O}{{\arabic*}})]
\item \label{obs1RM} Number of hospitalized people in ICU beds:
$$\observation(\state(\tau))=\hospitalizedC(\tau);$$

\item \label{obs2RM} Number of active cases (infectious persons detected):
$$\observation(\state(\tau))=\infected(\tau)+\hospitalized(\tau)+\hospitalizedC(\tau).$$
\end{enumerate}

As established in \cref{sec:triggered}, given an observation time window of $\hist =14$ days, the indicators to be used are defined for the recent history of the state
$$\stateh_{\hist}(t)=(\state(t-\hist),\ldots,\state(t)) \in \X^{\hist+1} .$$
Thus, the indicators considered in this example are computed from the observations defined in \ref{obs1RM} and \ref{obs2RM} as follows:
\begin{enumerate}[label=(\alph*)]
\item \label{indta} the mean of instantaneous observations:
$$\indicator(\state_{\hist}(t))= \frac{1}{\hist +1} \sum_{\tau=t-\hist}^t \observation(\state(\tau)) ;$$
\item  \label{indtb} the mean of differences of instantaneous observations:
$$\indicator(\state_{\hist}(t))= \frac{1}{\hist +1} \sum_{\tau=t-\hist+1}^t (\observation(\state(\tau)) -\observation(\state(\tau-1)))= \frac{\observation(\state(t))-\observation(\state(t-\hist))}{\hist +1} .$$
\end{enumerate}
Thus, we will assess four indicators associated with two instantaneous observations \ref{obs1RM} and \ref{obs2RM} and with two different ways to consider them \ref{indta} and \ref{indtb}.

Finally, to compare the performance of different event-triggered controls based on the indicators introduced previously (as introduced in \cref{sec:outcomes} and \cref{sec:comparison}), we consider two outcomes, that is, the function of outcomes  $\output: \X \times \UU \longrightarrow \R^2$ is given by $\output(\ic,\u)=(\outputc_1(\ic,\u),\outputc_2(\ic,\u))$, where the  event-triggered control $\u$ is defined by  $\u=\trigger(\hist,\delay,\indicator,\threshold)$,  for a given trigger threshold $\threshold$,  from an initial state $\ic$ at initial time $\it$. The outcomes considered are:
\begin{enumerate}[label=(\subscript{P}{{\arabic*}})]
\item \label{p1}Peak of ICU demand:  
$$ \outputc_1(\ic,\u)= \max_{\it\le t \le \horizon+1} \hospitalizedC(t) ;$$

\item  \label{p2} Total number of days in lockdown. This is expressed as a percentage (over the total simulation time) and is computed differently depending on whether our NPI is applied at the initial time $\it$ or not. If the NPI is applied  at $\it$, it is given by

$$\outputc_2(\ic,\u)=\frac{100}{(\horizon - \it) }\sum_{k=0}^{k_{\u}} (t_{2k+1}-t_{2k})$$
otherwise, we have
$$\outputc_2(\ic,\u)=\frac{100}{(\horizon - \it)}\sum_{k=0}^{k_{\u}} (t_{2k+2}-t_{2k+1}) .$$

We note that in the above expressions, $k_{\u}$ is the total number of the switches triggered  by the policy $\u$.
\end{enumerate}

\begin{figure}[ht]
\begin{center}
\scalebox{.72}{
\includegraphics{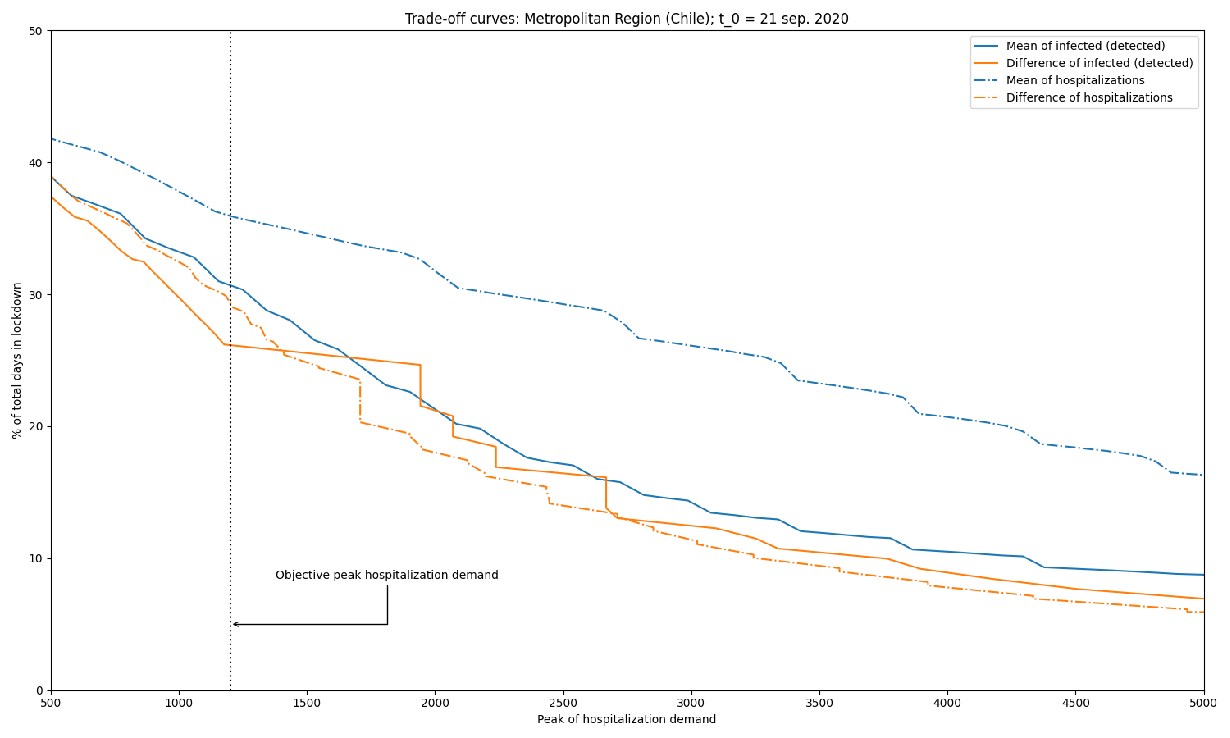}}
\end{center}
 \caption{Trade-off curves for our case study  of the spread of COVID-19 in Metropolitan Region, Chile, considering four indicators: number of hospitalized patients in ICU beds (observation \ref{obs1RM}) considering the mean (indicator \ref{indta}, blue dashed curve) and the mean of difference (indicator \ref{indtb}, orange dashed  curve) and the number of active cases (observation \ref{obs2RM}) considering the mean (indicator \ref{indta}, blue continuous curve) and the mean of difference (indicator \ref{indtb}, orange continuous curve).}\label{fig:rm}
\end{figure}

For a given event-triggered feedback  $\u=\trigger(\hist,\delay,\indicator,\threshold)$ we recall  the definition of the  trade-off curve associated with indicator $\indicator$, and parametrized by thresholds $\threshold$ (see \cref{sec:outcomes})
\begin{equation*}
\tradeoff(\ic,\hist,\delay,\indicator)=\{\output(\ic,\u) ~|~\u=\trigger(\hist,\delay,\indicator,\threshold);~~\threshold \in \thresholdd_{\indicator}\} \subseteq \R^2 ,
\end{equation*}
where $\thresholdd_{\indicator} \subseteq \R$ is the domain for the thresholds associated with indicator $\indicator$. For the four indicators considered  here (two instantaneous observations \ref{obs1RM} and \ref{obs2RM} and two different ways to consider them \ref{indta} and \ref{indtb}) and the two introduced outcomes \ref{p1} and \ref{p2}, the  trade-off curves are depicted in \cref{fig:rm}.

It is observed from  \cref{fig:rm} that the performance of the indicator corresponding to the mean of hospitalized patients in ICU beds (observation \ref{obs1RM} considered as  \ref{indta}), depicted by a blue dashed curve, is the worst in comparison to the other indicators in the sense that its trade-off curve is above the other curves; therefore, for any objective in the peak of ICU demand (outcome \ref{p1}), the use of this indicator will imply more days in the lockdown (outcome \ref{p2}). The trade-off curves of the other three indicators are quite similar. Nevertheless, the indicators based in the number of active cases (based in observation  \ref{obs2RM}, depicted by continuous curves) in practice can be less reliable that the indicator corresponding to the mean of the difference of the hospitalized patients in ICU beds (based in observation \ref{obs1RM} considered as  \ref{indtb}) that is shown by the orange dashed curve.

In \cref{table:rm},  we present the associated threshold for each indicator, having as objective a peak of ICU demand of $\outputc_1(\ic,\u)=1,200$ beds (see \ref{p1}). This value corresponds to an approximation of the total ICU beds available during October 2020 for the Metropolitan Region according to Chilean official data \cite{minciencia}. In this table, we also show the percentage of days in lockdown $\outputc_2(\ic,\u)$ (see \ref{p2}) associated with each indicators for the objective mentioned above.

\begin{table}[ht]
\vspace{.2cm}
\begin{center}
\footnotesize
\begin{tabular}{|l|c|c|}
\hlineB{3}
{\bf Indicator} & {\bf Threshold}  &  {\bf \% in lockdown}\\
  \hlineB{3}
$\indicator$  & $\threshold$  &  $\outputc_2(\ic,\u)$\\
\hlineB{3}
 Mean of ICU \ref{obs1RM}; \ref{indta} & 253& 36\%\cr
\hline
 Difference of ICU \ref{obs1RM}; \ref{indtb} & 0.4 & 29\%\cr
\hline
Mean of active cases \ref{obs2RM}; \ref{indta} & 87 & 31\%\cr
\hline
Difference of active cases \ref{obs2RM}; \ref{indtb} & 0.1 & 26\%\cr
\hlineB{3}
\end{tabular}
\vspace{3mm}\
\caption{Thresholds and percentages of days in lockdown $\outputc_2(\ic,\u)$ (see \ref{p2}) associated with four assessed indicators, considering  a peak of ICU demand objective of $\outputc_1(\ic,\u)=1,200$ beds (see \ref{p1}), in Metropolitan Region (Chile).}
\label{table:rm}
\end{center} 
\normalsize
\end{table}

\subsection{Case study: China}\label{sec:china}

In this example, we consider the compartmental model published in \cite{ivorra:2020} for the spread of COVID-19 in China.  The model is composed by nine state variables corresponding to different stages of the disease:   susceptible ($S$), exposed ($E$), infectious ($I$), infectious but undetected ($I_u$), hospitalized that will recover ($H_R$), hospitalized that will die ($H_D$), recovered after being detected ($R_d$), and recovered after being infectious but undetected ($R_u$), and dead by COVID-19 ($D$). Thus, the vector of the state variables is
$$\state(t) = (S(t),E(t),I(t), I_u(t), H_R(t), H_D(t), R_d(t), R_u(t),D(t)) \in \X=\R^9_+ .$$

Since the model in \cite{ivorra:2020} is established in continuous-time, we consider a very small time step (one hour). Then, we use  the standard Euler method to discretize the  dynamics in \cite{ivorra:2020}  with a step $\paso=1/24$ to obtain $\dynamics:\X\times\U\longrightarrow\X$ as in \eqref{eq:system}.
Thus, for this case-study, the (autonomous) dynamics that define the discrete-time control systems are described by the following: 

\begin{equation}
\label{eq:modelf_China}
\dynamics(\state,\control)=\dynamics(\susceptible,\exposed,\infected, \infected_u ,\hospitalized_R, \hospitalized_D, \recovered_d,\recovered_u,\dead,\control) =
\left(\begin{array}{c}

\susceptible -\paso \hat \rateg(\state,\control)\susceptible 
 \\[3mm]

 \exposed+ \paso \left( \hat \rateg(\state,\control)\susceptible - \duration_{\exposed}\exposed \right) \\[3mm]

\infected+  \paso \left( \duration_{\exposed} \exposed -\duration_{\infected}\infected   \right) \\[3mm]

 \infected_u +   \paso \left( (1-\fraction_{\infected\hospitalized_R}-\fraction_{\infected\hospitalized_D})\duration_{\infected} \infected -\duration_{\infected_u}\infected_u   \right)  \\[3mm]

 \hospitalized_R +  \paso \left( \fraction_{\infected\hospitalized_R} \duration_{\infected} \infected - \duration_{\hospitalized_R}\hospitalized_R \right) \\[3mm]

 \hospitalized_D + \paso \left( \fraction_{\infected\hospitalized_D} \duration_{\infected} \infected
 - \duration_{\hospitalized_D} \hospitalized_D  \right) \\[3mm]

 \recovered_d+  \paso \duration_{\hospitalized_R} \hospitalized_R \\[3mm]

 \recovered_u +  \paso \duration_{\infected_u} \infected_u \\[3mm]

 \dead+ \paso  \duration_{\hospitalized_D} \hospitalized_D 
\end{array}\right) .
\end{equation}

Here, the controlled contagion rate $\hat \rateg(\state,\control)$  is given by 
\begin{equation}\label{eq:controlledrate-china}
\hat \rateg(\state(t),\control(t))=\frac{1}{N} (1- \control(t))
(\rate_{\exposed} \exposed(t)  +\rate_{\infected} \infected(t) + \rate_{\infected_u} \infected_u(t) + \rate_{\hospitalized_R}\hospitalized_R(t)+\rate_{\hospitalized_D}\hospitalized_D(t))  ,
\end{equation}
where a single control variable  $\control \in \U=[0,1-\hat \delta]$ has been considered. This control is associated with the implementation of NPIs such as lockdowns. It multiplies all of the contagion rates ($\rate_{\exposed}$, $ \rate_{\infected}$, $ \rate_{\infected_u}$, $\rate_{\hospitalized_R}$ and $\rate_{\hospitalized_D}$), because its main effect is to reduce the contact rates among all individuals. Parameter $\hat \delta \in (0,1)$ represents the fraction of the population that reduces their contact rates during a lockdown (or more generally the application of a given NPI).  The latter is considered to be different from zero because, even during lockdowns,  some basic services must still operate. In the following, we set $\hat \delta=0.25$.

The parameters defining the dynamics $\dynamics:\X\times\U\longrightarrow\X$ in \eqref{eq:modelf_China} are given by the vector $(\rate,\duration,\fraction) \in \R^5_+ \times \RR^5_+  \times [0,1]^2$. Vector $\rate = (\rate_{\exposed}, \rate_{\infected}, \rate_{\infected_u},\rate_{\hospitalized_R}, \rate_{\hospitalized_D}) \in \R^5_+$ contains all of the specific rates associated with the contagious stages of the disease. Parameters $\duration=(\duration_{\exposed}, \duration_{\infected}, \duration_{\infected_u},\duration_{\hospitalized_R},\duration_{\hospitalized_D}) \in \R^5_+$ are the mean rates of the transition from these respective stages to the subsequent stage. In other words, for a disease stage $X \in \{\exposed,\infected,\infected_u,\hospitalized_R,\hospitalized_D\}$ parameter $\duration_X^{-1}$ days represents the mean duration of stage $X$.  Finally, the vector $\fraction=(\fraction_{ \infected \hospitalized_R},\fraction_{\infected \hospitalized_D}) \in [0,1]^2 $ allows us to describe  the distribution from infected people ($\infected$) to the next three stages: infected but undetected  ($\infected_u$),  hospitalized that will recover ($H_R$), and hospitalized that will die ($H_D$). Indeed, the fraction of infected people that are not detected is given by $1-\fraction_{ \infected \hospitalized_R}-\fraction_{\infected \hospitalized_D}$, while $\fraction_{ \infected \hospitalized_R}$ and $\fraction_{ \infected \hospitalized_D}$ represent the fraction passing from infected to the hospitalized stages $H_R$ and $H_D$, respectively.  We note that in \cite{ivorra:2020}, both fractions are equivalently written in terms of the case fatality ratio and of the fraction of infected people that are detected. The values of $(\rate,\duration,\fraction)$ used in our simulations correspond to those used in  \cite{ivorra:2020} for experiment  {\bf EXP$_{\rm \bf 29M}$}; see Table 3 in \cite{ivorra:2020}.

Additionally, the initial conditions 
\begin{align*}
\state(\it) &= (\susceptible_0,\exposed_0,\infected_0, {\infected_u}_0 ,{\hospitalized_R}_0, {\hospitalized_D}_0, {\recovered_d}_0,{\recovered_u}_0,\dead_0)
\in \X = \R^9_+ 
\end{align*}
were obtained  from the calibration realized with the complete set of data available in  \cite{ivorra:2020}, that is,  our initial condition vector $\state(\it)$ corresponds to the last day output (March 29, 2020) of the fitted model in  \cite{ivorra:2020} for experiment  {\bf EXP$_{\rm \bf 29M}$}. These values are given in \cref{table:initial_conditions_china}.

\begin{table}[ht]
\vspace{.2cm}
\begin{center}
\footnotesize
\begin{tabular}{|l|r|}
\hlineB{3}
{\bf State variable} & {\bf Value} \\
\hlineB{3}
 $\susceptible_0$ & 1,389,828,000  \cr
 \hline
 $\exposed_0$ & 14 \cr
 \hline
 $\infected_0$ & 2 \cr
 \hline
$  {\infected_u}_0$& 1,555 \cr
 \hline
$ {\hospitalized_R}_0$&  2,035  \cr
 \hline
$ {\hospitalized_D}_0$& 270 \cr
 \hline
$ {\recovered_d}_0$& 73,622 \cr
 \hline
$ {\recovered_u}_0$&  90,346 \cr
\hline
$\dead_0$&  3,708 \cr
\hlineB{3}
\end{tabular}
\vspace{3mm}\
\caption{Initial conditions for example in \cref{sec:china} corresponding to China, estimated at 
$\it=$ March 29, 2020.
}
\label{table:initial_conditions_china}
\end{center} 
\normalsize
\end{table}

It is important  to note that a very strict lockdown was considered in the experiment  {\bf EXP$_{\rm \bf 29M}$} in   \cite{ivorra:2020}. This explains the very low values obtained for the infected stages in \cref{table:initial_conditions_china}. Recall that our initial time ($\it=$ March 29, 2020) corresponds to the final time in that experiment after the application of the lockdowns.

Regarding the other elements necessary to pursue our analysis,  the event-triggered controls introduced in \cref{sec:triggered} are defined similarly to those  in  \cref{sec:chile}. Indeed, the considered \emph{observation time window} and the \emph{minimal time of implementation} (or minimal interevent time) will be $\hist=\delay=14$ days. The event-triggering mechanism described in \eqref{eq:controller0} and  \eqref{eq:controllerk} is defined by $\control_{\rm ref}=1-\hat \delta$ (associated with $\control(\it)$ in \eqref{eq:controller0})  and  by \emph{controllers} $\hat \control_-,~\hat \control_+: \titf{0}{\horizon-\it} \times \U \longrightarrow  \U$ (see  \eqref{eq:controllerk}) that are considered as max-linear functions  of the form
\begin{eqnarray*}
\hat \control_-(\tau,\control)&=&\max\left\{0, \control\left(1 - \frac{\tau}{\delay}\right)\right\}\\[2mm]
\hat \control_+(\tau,\control)&=&\min\left\{1 -\hat \delta, \control\left(1 - \frac{\tau}{\delay}\right)  + \left(\frac{1-\hat \delta}{\delay}\right) \tau\right\} .
\end{eqnarray*}
As in the previous case-study analyzed in \cref{sec:chile}, our control saturates in $\delay$ days, coinciding with the minimum interevent time. Thus, we assume that at the end of  this period of time, the effect of NPI (or their release) is fully accomplished, producing maximal isolation ($\control = 1-\hat \delta$) or reaching the normal contact rates equal to those prior to the pandemic ($\control = 0$). 

To define trigger indicators $\indicator: \X^{\hist+1} \longrightarrow \R$ and then the event-triggered sets $\Set$ (cf. \cref{sec:tradeoff}), we also consider the following instantaneous observations of the state:
\begin{enumerate}[label=(\subscript{\tilde O}{{\arabic*}})]
\item \label{obs1:china} Number of hospitalized people:
$$\observation(\state(\tau))=H_R(\tau)+H_D(\tau);$$

\item \label{obs2:china} 
Number of detected infectious people:
$$\observation(\state(\tau))=I(\tau).$$
\end{enumerate}

As established in \cref{sec:triggered}, given an observation time window of 14 days (i.e., $\hist=14$), the indicators to be used are defined for the recent history of the state
$$\stateh_{\hist}(t)=(\state(t-\hist),\ldots,\state(t)) \in \X^{\hist+1} .$$
Thus, the indicators considered in this case-study are computed from observations defined in \ref{obs1:china} and \ref{obs2:china} as follows:

\begin{itemize}
\item[(a)] the mean of instantaneous observations:
$$\indicator(\state_{\hist}(t))= \frac{1}{\hist +1} \sum_{\tau=t-\hist}^t \observation(\state(\tau)) ;$$
\item[(b)] the mean of differences of instantaneous observations:
$$\indicator(\state_{\hist}(t))= \frac{1}{\hist +1} \sum_{\tau=t-\hist+1}^t (\observation(\state(\tau)) -\observation(\state(\tau-1)))= \frac{\observation(\state(t))-\observation(\state(t-\hist))}{\hist +1} .$$
\end{itemize}

Finally, for the assessment of the different event-triggered controls policies based on the indicators above, we apply the methodology described in  \cref{sec:outcomes} and \cref{sec:comparison}. Thus, we consider the function of outcomes  $\output: \X \times \UU \longrightarrow \R^2$ given by $\output(\ic,\u)=(\tilde \outputc_1(\ic,\u),\tilde \outputc_2(\ic,\u))$, where the  event-triggered control $\u$ is defined by  $\u=\trigger(\hist,\delay,\indicator,\threshold)$ (for a given trigger threshold $\threshold$,  from an initial state $\ic$ at initial time $\it$), and the outcomes $\tilde \outputc_1(\ic,\u)$ and $\tilde \outputc_2(\ic,\u)$ are given by the following:

\begin{enumerate}[label=(\subscript{\tilde P}{{\arabic*}})]
\item \label{p1:china} Peak of hospitalization demand: 
$$ \tilde \outputc_1(\ic,\u)= \max_{\it\le t \le \horizon+1}(H_R^{\u}(t)+H_D^{\u}(t)) ,$$
where $H_R^{\u}(t) + H_D^{\u}(t)$ corresponds to the total number of hospitalized patients at period $t$.

\item \label{p2:china}  Total number of days of lockdown: This outcome was already defined in \cref{sec:chile}, outcome \ref{p2}. For the sake of completeness, we recall its computation here:

If the NPI is applied  at $\it$, it reads as follows
$$\tilde \outputc_2(\ic,\u)=\frac{100}{(\horizon - \it) }\sum_{k=0}^{k_{\u}} (t_{2k+1}-t_{2k})$$
otherwise, we have
$$\tilde \outputc_2(\ic,\u)=\frac{100}{(\horizon - \it)}\sum_{k=0}^{k_{\u}} (t_{2k+2}-t_{2k+1}) .$$
We note that in the previous expressions, $k_{\u}$ is the total number of the switches triggered  by the policy $\u$.
\end{enumerate}

\begin{figure}[ht]
\begin{center}
\scalebox{.72}{
\includegraphics{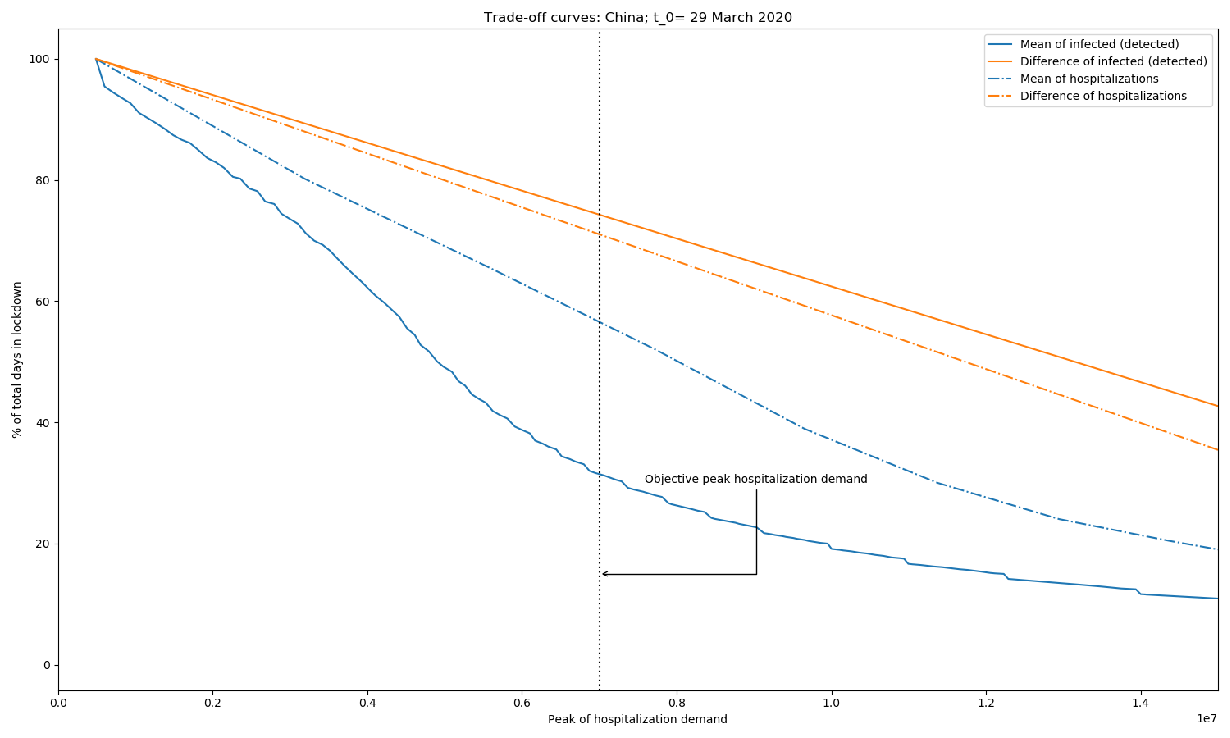}}
\end{center}
 \caption{Trade-off curve for our case-study based on \cite{ivorra:2020} for the spread of COVID-19 in China. Four indicators are considered: nmber of hospitalized people (observation \ref{obs1:china}) considering the mean (indicator \ref{indta}, blue dashed curve) and the mean of difference (indicator \ref{indtb}, orange dashed  curve) and the number of active cases (observation \ref{obs2:china}) considering the mean (indicator \ref{indta}, blue continuous curve) and the mean of difference (indicator \ref{indtb}, orange continuous curve).
 }\label{fig:curve-China}
\end{figure}

Hence, the trade-off curves (see \cref{sec:outcomes})
\begin{equation*}
\tradeoff(\ic,\hist,\delay,\indicator)=\{\output(\ic,\u) ~|~\u=\trigger(\hist,\delay,\indicator,\threshold);~~\threshold \in \thresholdd_{\indicator}\} \subseteq \R^2 ,
\end{equation*}
for the four indicators introduced above  (two instantaneous observations \ref{obs1:china} and \ref{obs2:china} and to the two different ways to consider them \ref{indta} and \ref{indtb}) are depicted in  \cref{fig:curve-China}.

\cref{fig:curve-China} shows that the 14-day rolling average of infected and detected people (observation \ref{obs2:china} considered as  \ref{indta}) is clearly the best indicator among the four analyzed indicators. The second best indicator is the mean of total hospitalizations due to COVID-19 (observation \ref{obs1:china} considered as  \ref{indta}), but its performance is considerably weaker. For instance, when we consider an objective peak hospitalization demand of 7 million of beds (that is, $\tilde \outputc_1(\ic,\u)=7,000,000$; see outcome \ref{p1:china}), the percentage of time in lockdowns (outcome \ref{p2:china}) for the policies associated with each indicator is given by $\tilde \outputc_2(\ic,\u)= $ 31\% and 57\%, respectively, over the total analyzed time. This demonstrates the strong advantage of using  a policy associated with the rolling average of infected people with respect to the other analyzed indicators in this case. Therefore, our analysis suggest the use of this indicator as the basis for an event-triggered policy for this case study.

In Table \cref{table:China} we show the percentage of the days in lockdown $\tilde \outputc_2(\ic,\u)$ (see \ref{p2:china}) associated with each indicator when the objective peak hospitalization demand is $\tilde \outputc_1(\ic,\u)=7,000,000$ beds (see \ref{p1:china}).
We also present the associated thresholds $\threshold$ for these four indicators.

\begin{table}[ht]
\vspace{.2cm}
\begin{center}
\footnotesize
\begin{tabular}{|l|c|c|}
\hlineB{3}
{\bf Indicator} & {\bf Threshold}  &  {\bf \% in lockdown}\\
  \hlineB{3}
$\indicator$  & $\threshold$  &  $\tilde \outputc_2(\ic,\u)$\\
\hlineB{3}
 Mean of hospitalized people \ref{obs1:china}; \ref{indta} &  21 & 57\%\cr
\hline
Difference of hospitalized people \ref{obs1:china}; \ref{indtb} & -1.29 & 71\%\cr
\hline
Mean of detected infectious people  \ref{obs2:china}; \ref{indta} & 3,935,670 & 31 \%\cr
\hline
 Difference of detected infectious people \ref{obs2:china}; \ref{indtb} & -615 & 74 \%\cr
\hlineB{3}
\end{tabular}
\vspace{3mm}\
\caption{Thresholds and percentages of days in lockdown $\tilde \outputc_2(\ic,\u)$ (see \ref{p2:china}) associated with the four assessed indicators considering  the peak hospitalization demand of $\tilde \outputc_1(\ic,\u)=7,000,000$ beds (see \ref{p1:china}) in China as the objective.}
\label{table:China}
\end{center} 
\normalsize
\end{table}

\section{Final remarks}
\label{sec:conclusions}

In this work, we addressed two main goals: to propose a modeling framework for representing the decision-making process related to the application of NPI policies based on the observation of epidemiological indicators and to introduce a methodology to compare the effectiveness of these strategies (defined by trigger indicators) with respect to multiple objectives outcomes. The proposed framework is that of event-triggered controls with some specific details such  as considering the state in a window of time and imposing a minimum intervention time (or minimum interevent time). Both features were introduced in order to represent the application of the NPIs more realistically.

The two case-studies discussed in \cref{sec:examples} show how our methodology can be used to recommend an indicator-based policy when some given outcomes are treated as objectives in our analysis. Nevertheless, our approach should be consider only as a guide to suggest a recommendation. In practice, decision makers must also undertake a complementary and more in-depth analysis. For instance, the observability of the indicators under analysis is something that must also be taken into account. Indeed, when an indicator based on active cases has a similar performance to that of another indicator based, for instance, on the number of hospitalized patients (as in the case-study given in \cref{sec:chile}), a possible recommendation is to consider a policy based on the second type of indicator because the corresponding observations (hospitalized patients) are more reliable and, consequently, the latter is a more robust NPI policy in practice.

As described in this article, our approach permits to represent the decision-making process based on the observation of indicators and compare the use of these indicators to define NPI policies from a cost-benefit perspective. This comparison will depend on the outcomes chosen as objectives and also on the structure of the dynamic model and on the indicators chosen for this comparison. Both of the latter mathematical objects are intrinsically related in the sense that the dynamics considered  as a basis in our methodology must allow the computation and simulation of the indicators that we seek to compare in the analysis. Therefore, for assessment of more elaborate indicators such as positivity rates, indicators associated with testing and contact tracing and effective reproduction number, it is necessary to consider dynamical systems with compartmental stages that allows the formulation of such indicators. On the other hand, our methodology was demonstrated  for two SEIR-type models that  are some of the most common mathematical models in epidemiology. Nevertheless, our methodology can also be adapted to other underlying mathematical models, such as agent-based models. Of course, different dynamics or other underlying mathematical models can lead to different recommendations. 
However, this issue is inherent to any real-life decision derived from a mathematical analysis.

The framework introduced in this paper for assessing indicators-based event-triggered policies raises some interesting problems to study the  further theoretical analysis of which can have important repercussions on the decision-making process. We mention three such problems:
\begin{itemize}

\item {\bf Time consistency:} When the analysis introduced in \cref{sec:comparison} for determining the best indicator to use is carried out at time $\it$, following the policy induced by this indicator, there is no guarantee that this strategy will be the best in the long term. Since the trade-off curves are dynamic objects because they depend on the initial conditions (see \eqref{eq:curve}), it is possible that with time,  another indicator should be chosen for triggering NPIs, because its performance is better considering the future state as the initial condition.

 In the terms of trade-off curves and the comparison process proposed in \cref{sec:comparison}, the aforementioned situation occurs when two of these curves (the best chosen at $\it$ and other), viewed as curves that evolve in time, intersect in a future time at the level (in the outcomes space) corresponding to the $m-1$ outcomes objectives for which the analysis was carried out at time $\it$. To prove that a future intersection between the trade-off curves will not occur will clearly require us to assume certain properties of the dynamics and the compared indicators; this is an interesting problem. A similar problem in the context of linear system where the performance (i.e., our outcomes) is a quadratic cost, is studied in \cite{Antunes:2016}.
 
 On the other hand, it is natural to consider situations where the outcomes objectives can evolve in time. For instance, at time $\it$ the situation is analyzed considering a given ICU maximal capacity but perhaps in some months, this capacity will be increased. In summary, the comparison process of indicators-based event-triggered policies should be implemented continuously across the evolution of the disease in the population, taking into account eventual updating of objective outcomes by the decision makers. Due the latter and the dynamics aspect of trade-off curves, with longer time it is possible that a discarded indicator will in fact become the best indicator.
 
 \item {\bf Is it better to consider policies based on more indicators?} In the analysis proposed in this work, we have considered only policies based on a single indicator. That is, to consider event-triggered sets $\Set$, defined in the space of the recent history of states $\X^{\hist+1}$, in the form
\begin{equation}\label{eq:setS}
\Set=\Set(\hist,\indicator,\threshold)=\{\stateh_{\hist} \in \X^{\hist+1}~|~\indicator(\stateh_{\hist}) \le \threshold\} ,
\end{equation}
where the indicator $\indicator: \X^{\hist+1} \longrightarrow \R$ takes values in $\R$ and $\threshold \in \R$ is a threshold for this indicator. Consideration of more indicators leads to the consideration of a vector-valued function $\indicator: \X^{\hist+1} \longrightarrow \R^p$, with $p \ge 2$, and the inequality $\indicator(\stateh_{\hist}) \le \threshold$  in \eqref{eq:setS} with respect to the componentwise order for a vector of thresholds $\threshold \in \R^p$. We wonder if the question of whether more indicators lead to better policies can be formulate as follows: If $\Set(\indicator, \threshold) \subseteq \Set(\tilde \indicator, \threshold)$ for all vector of threshold $\threshold$ (in an appropriate domain) can one conclude that a policy based on $\indicator$ is better than one based on $\tilde \indicator$ (or vice versa)?

For triggering NPIs it can be natural to consider various indicators, as in the examples cited in the introductory section, motivated by the reliability (or the lack of) of observations for making decisions. Nevertheless, assuming that all observations are reliables, it is not clear that the addition of indicators for triggering the decisions will induce better policies  from the cost-effectiveness viewpoint.

\item {\bf Stochastic analysis:} The framework introduced in this work can be extended considering stochastic dynamics and indicators. In this case, a possible additional outcome that can be considered is the maximum  probability of satisfying the other outcomes (seen as objectives), and then the performance of different event-triggered policies can be assessed in terms of that probability, considering the risk-aversion of the decision makers.

To considering stochastic dynamics,  the estimated probabilities distributions of involved parameters can be take as  one of the output given by the calibration procedure used in \cref{sec:chile} and mentioned in \cref{app:chile}.
\end{itemize}

Perhaps using simple models (few state variables) it is possible to provide mathematical proofs for the questions posed in the first two points of the previous list, that is, to ensure that the trade-off curves do not intersect and to deduce when an indicator is better than others from the properties of the event-triggered sets $\Set$ corresponding to those indicators. Eventually the theoretical analysis can be simpler considering feedbacks based only on the instantaneous observations and time-continuous dynamics, a framework have purposefully avoided here because the measurability of the involved functions merits a deeper analysis. The aforementioned problems are beyond the scope of this paper; however, we plan to study them in future works.

\appendix
\section{Parameter description and calibration procedure of model in  \cref{sec:chile}}\label{app:chile} 

The parameters of system \eqref{eq:model} associated with dynamics $\dynamics:\X\times\U\longrightarrow\X$ given by \eqref{eq:modelf} are contained in the vector
$(\rate,\duration,\fraction) \in \R^3_+ \times \RR^5_+  \times [0,1]^5$, where $\rate = (\rate_{\exposed}, \rate_{\infectedU}, \rate_{\infected}) \in \R^3_+$ is the vector of the specific rates associated with contagious stages of the disease and that are involved in the definition of the controlled contagion rate $\hat \rateg(\state(t),\control(t))$ given by \eqref{eq:controlledrate}. 

Parameters $\duration=(\duration_{\exposed}, \duration_{\infectedU}, \duration_{\infected},\duration_{\hospitalized},\duration_{\hospitalizedC}) \in \R^5_+$ are the mean rates of transition from the respective stages to the next stage. In other words, for a disease stage $X \in \{\exposed,\infectedU,\infected,\hospitalized,\hospitalizedC\}$ parameter $\duration_X^{-1}$ days represents the mean duration of stage $X$. 

The vector $\fraction=(\fraction_{\exposed \infected},\fraction_{\infected \recovered},\fraction_{\hospitalized \recovered},\fraction_{\hospitalized\dead},\fraction_{\hospitalizedC\dead}) \in [0,1]^5 $ contains different distribution fractions between the different stages. Thus, $\fraction_{\exposed \infected}$ is the fraction of exposed people who become infected (with symptoms), $\fraction_{\infected \recovered}$ is the fraction of infected people who recover, $\fraction_{\hospitalized \recovered}$ is the fraction of hospitalized (in normal services) people who recover, $\fraction_{\hospitalized\dead} $ is the fraction of hospitalized (in normal services) people who die, and  $\fraction_{\hospitalizedC\dead} $ is the fraction of hospitalized people in ICU beds who die.

The parameters  have been estimated using a Hamiltonian Markov chain Monte Carlo estimation method using the Stan library \cite{Carpenter:2017}. The method is a Bayesian approach that allows the approximation of the probability densities for each of the parameters in the system when using a sufficient number of iterations. We use the Gelman-Rubin ($\hat{R}$) convergence diagnostic and 12 parallel chains to ensure that the separate chains converge to the same values. $\hat{R}$ compares the variance of each chain to the pooled variance of all chains. Each of the parameters obtained an $\hat{R} \leq 1.03$ where each chain ran 10,000 iterations of which a half were warm-up.

For each of the parameters, we set upper and lower limits to restrict the parameter search space to allow realistic combinations only. In addition, we used a weakly informative prior for each of the parameters with mean values obtained from the literature and broad standard deviations.

 In \cref{table:parameters}, we state for each parameter its limits, prior probability distribution in the form of a mean and standard deviation, and some references for the mean of the used prior distributions. The model was trained on data  from June 20, 2020, until September 20, 2020, in the Metropolitan Region (Chile), considering  the total number of COVID-19 cases (symptomatic), the number of patients in ICU, and the number of deceased, according the information collected and published by the Chilean Ministry of Science, Technology, Knowledge and Innovation in \cite{minciencia}. We assume an  independent and identically distributed (i.i.d.) error of $\sigma$ for all data points.
 
 Since there are no data available of the asymptomatic population, we are not able to train the trajectory of asymptomatic cases from $\exposed$ to $\infectedU$ and to $\recovered$. That is, the parameters $\fraction_{\exposed \infected}$ and $\duration_{\infectedU}$ are not trained by maximizing the posterior, except for the small influence through $\rate_{\infectedU}$ that is difficult to distinguish from the other sources of contagion. To mitigate this, we will fix $\fraction_{\exposed \infected} = 0.6$ where we assume 60\% of the infections to be symptomatic,  following \cite{Oran:2020} and \cite[Scenario 5: Current Best Estimate]{cdcbest}.
 
\begin{table}[ht]
\begin{center}
\footnotesize
    \begin{tabular}{|c|c|c|c|c|}
\hlineB{3}
    {\bf Parameter} & {\bf Limits} & {\bf Prior} & {\bf Value} & {\bf Reference (mean prior)}\\[1mm]
  \hlineB{3}
$\rate_{\exposed}$ & $[0.002, 0.2]$ & $\mathcal{N}(0.016, 0.1)$ & $0.04\pm\num{3.3e-5}$ &----- \\[1mm]
 \hline
$\rate_{\infectedU}$ & $[0.002, 0.2]$ & $\mathcal{N}(0.016, 0.1)$ & $0.04\pm\num{3.3e-5}$ &-----  \\[1mm]
 \hline
$\rate_{\infected}$ & $[0.01, 1]$ & $\mathcal{N}(0.08, 0.5)$ & $0.2\pm\num{1.7e-4}$&----- \\[1mm]
 \hline
$\duration_{\exposed}$ & $[1/15, 1/2]$ & $\mathcal{N}(1/5.8,0.07 )$ & $0.39\pm\num{3.9e-4}$& \cite{McAloone:2020}  \\[1mm]
 \hline
$\duration_{\infectedU}$ & $[1/30, 1/5]$ & $\mathcal{N}(1/9.5, 0.07)$ & $0.17\pm\num{1.2e-4}$& \cite{Byrnee:2020}  \\[1mm]
 \hline
$\duration_{\infected}$ & $[1/30, 1/5]$ & $\mathcal{N}(1/7,0.07 )$ & $0.17\pm\num{1.4e-4}$& \cite{Wang:2020}  \\[1mm]
 \hline
$\duration_{\hospitalized}$ & $[1/30, 1/5]$ & $\mathcal{N}(1/11, 0.07)$ & $0.17\pm\num{1.1e-4}$ & \cite{Zhou:2020} \\[1mm]
 \hline
$\duration_{\hospitalizedC}$ & $[1/30, 1/5]$ & $\mathcal{N}(1/8,0.07 )$ & $0.14\pm\num{1.1e-4}$& \cite{Zhou:2020}  \\[1mm]
 \hline
$\fraction_{\infected\recovered}$ & $[0.5, 0.99]$ & $\mathcal{N}(0.7395, 0.5 )$ & $0.61\pm\num{3.2e-4}$& \cite{covid2020severe} \\[1mm]
 \hline
$\fraction_{\hospitalized\recovered}$ & $[0.5, 0.99]$ & $\mathcal{N}(0.6585, 0.5)$ & $0.61\pm\num{3e-4}$&\cite{Wang:2020}  \\[1mm]
 \hline
$\fraction_{\hospitalized\dead}$  & $[0.01, \min\{0.5, 1-\phi_{HR}\}]$ & $\mathcal{N}(0.0244, 0.5 )$ & $0.12\pm\num{1.5e-4}$& \cite{Wang:2020} \\[1mm]
 \hline
$\fraction_{\hospitalizedC\dead}$ & $[0.01, 0.5]$ & $\mathcal{N}(0.385, 0.5)$ & $0.12\pm\num{4.1e-4}$& \cite{Huang:2020} \\[1mm]
 \hline
$\sigma$ & $[0.0, 1.0]$ & ${\rm Cauchy}(0, 1)$ & $0.02\pm\num{3e-6}$& ----- \\[1mm]
  \hlineB{3}
    \end{tabular}
    \caption{Parameter limits, priors, and mean and standard deviation following MCMC estimation.}
    \label{table:parameters}
\end{center} 
\normalsize
\end{table}

The initial conditions $\state(\it)=(\susceptible_0,\exposed_0,\infectedU_0,\infected_0,\recovered_0,\hospitalized_0,\hospitalizedC_0,\dead_0) \in \X = \R^8_+$, 
have been estimated using weakly informative priors where we take $\exposed_0$, $\infectedU_0$, $\infected_0$, and $\hospitalized_0$ to have scaled Beta priors between $0$ and $N$ with $\alpha=1$ and $\beta=4$ (parameters of Beta distribution), where $N = 7,112,808$ is the total population. For the initial number of hospitalized  patients in ICU, deceased patients, and total number of infections to date we assume a lognormal distribution with a standard deviation of $\sigma$ around the data. In this setup, $\susceptible_0$ and $\recovered_0$ are free variable but we need to restrict at least one of these variables  to solve the system. We propose a rough estimate of $\recovered_0$ to be the total number of infected patients to date plus the percentage that were mildly infected. Here, we assume that the number of people that are currently infected, hospitalized, or deceased are small relative to the total population. $\susceptible_0$ then follows from $\susceptible_0 = N - \exposed_0 - \infectedU_0 - \infected_0 - \hospitalized_0 - \hospitalizedC_0 - \recovered_0 - \dead_0$. The obtained values of initial conditions are shown in \cref{table:ic}.

\begin{table}[!ht]
\begin{center}
\footnotesize
\begin{tabular}{|c|c|}
\hlineB{3}
\bf State &\bf  Initial condition\\
\hlineB{3}
$\susceptible_0$&  6,671,557\\
\hline
$\exposed_0$& 1,697  \\
\hline
$\infectedU_0$& 1,723  \\
\hline
$\infected_0$&  2,540\\
\hline
$\hospitalized_0$& 1,157 \\
\hline
$\hospitalizedC_0$& 433 \\
\hline
$\recovered_0$&421,948  \\
\hline
$\dead_0$& 11,753 \\
\hlineB{3}
\end{tabular}
\vspace{3mm}\
\caption{Initial conditions for  example in \cref{sec:chile} corresponding to Metropolitan Region (Chile), estimated at $t_0=$ September 21, 2020.}\label{table:ic}
\end{center}
\normalsize
\end{table}

\section*{Acknowledgments}
We are very grateful to Benjamin Ivorra  (Complutense University of Madrid, Spain) for his enlightening advice and assistance for  the example introduced in \cref{sec:china}. Any flaw or mistake in this section is our responsibility. 

\bibliographystyle{siamplain}
\bibliography{references}

\end{document}